\documentclass[11 pt, twoside]{amsart}
\usepackage{parskip}
\usepackage{amsthm,amsmath,amssymb,oldgerm}
\usepackage[a4paper, margin=2.5cm]{geometry}
\usepackage{mathrsfs}
\usepackage{verbatim}
\usepackage{hyperref}
\theoremstyle{plain}
\newtheorem{thm}{Theorem}[section]
\newtheorem*{theorem*}{Main theorem}
\newtheorem{mainthm}[]{Main Theorem}[]
\newtheorem{lem}[thm]{Lemma}
\newtheorem{prop}[thm]{Proposition}

\theoremstyle{definition}
\theoremstyle{definition}
\newtheorem{rem}[thm]{Remark}
\def\Bn{\mathbb{B}^{n}}

\def \Q{\mathbb{Q}}

\def \G {\Gamma}

\def \C {\mathbb{C}}

\def \R {\mathbb{R}}

\def \D {\mathbb{D}}

\def\E{\mathbb E}

\def\Z{\mathbb Z}

\def \caln {\mathcal{N}}
\def \calf {\mathcal{F}}

\def \calh {\mathcal{H}}

\def \cala {\mathcal{A}}

\def \calo {\mathcal{O}}
\def \cals {\mathcal{S}}

\DeclareMathOperator{\bkn} {{\it{\mathcal{B}_{\Gamma}^{k(n+1)}}}}

\DeclareMathOperator{\bkx} {{\it{\mathcal{B}}}_{{\it{X}}_{\G}}^{\lambda^{{\it k}}}}

\DeclareMathOperator{\berkvol}{\mu_{\mathrm{Ber},{\it{k}}}^{\mathrm{vol}}}

\DeclareMathOperator{\ckn}{{\it{C_{\mathrm{X_{\Gamma}}}^{k(n+1)}}}}
\DeclareMathOperator{\dhyp}{d_{\mathrm{hyp}}}
\DeclareMathOperator{\hyp}{\mu_{\mathrm{hyp}}}

\DeclareMathOperator{\hypbn}{\mu_{hyp}}
\DeclareMathOperator{\hypbnvol}{\mu_{hyp}^{vol}}

\DeclareMathOperator{\sk}{{\it{S_{k}}}(\Gamma)}
\DeclareMathOperator{\skn}{{\it{S_{k(n+1)}}}(\Gamma)}

\DeclareMathOperator{\rx}{{\it{r_{X_{\mathrm{\G}}}}}}
\title[Estimates of Picard automorphic forms]{Estimates of automorphic forms on $\mathrm{SU}(n,1)$}
\author{Anilatmaja Aryasomayajula}
\address{Department of Mathematics, Indian Institute of Science Education and Research (IISER) Tirupati, 
Transit campus at Sri Rama Engineering College, Karkambadi Road,
Mangalam (B.O),Tirupati-517507, India.}
\email{anil.arya@iisertirupati.ac.in}
\author{Baskar Balasubramanyam}
\address{Department of Mathematics, Indian Institute of Science Education and Research (IISER) Pune,
Dr. Homi Bhabha Road, Pashan, Pune 411008, India.}
\email{baskar@iiserpune.ac.in}
\date{\today}
\begin{document}
%%%%%%%%%%%%%%%%%%%%%%%%%%%%%%%%%%%%%%%%%%%%%%%%%%%%%%%%%%%%%%%%%%%%%%%%%%%%%%%%%%%%%%%%%%%
\begin{abstract}
For a fixed $n\geq 2$, let $\G\subset \mathrm{SU}((n,1),\calo_{K})$ be a torsion-free, finite-index subgroup, where $\calo_K$ denotes the ring of integers of a totally imaginary number field $K$ of degree $2$. We assume that $\G$ admits only one cusp at $\infty$. Let $\Bn$ denote the $n$-dimensional complex ball endowed with the hyperbolic metric 
$\mu_{\mathrm{hyp}}$, and let $X_{\G}:=\G\backslash \Bn$ denote the quotient space, which is a noncompact complex manifold of dimension $n$. Furthermore, let 
$\mu_{\mathrm{hyp}}^{\mathrm{vol}}$ denote the volume form associated to the $(1,1)$-form $\mu_{\mathrm{hyp}}$.
%%%%%%%%%%%%%%%%%%%%%%%%%%%%%%%%%%%%%%%%%%%%%%%%%%%%%%%%%%%%%%%%%%%%%%%%%%%%%%%%%%%%%%%%%%%
 
\vspace{0.1cm}\noindent
Let $\Lambda:=\Omega_{\overline{X}_{\G}}^{n}$ denote the line bundle, whose sections are holomorphic $(n,0)$-forms, where $\overline{X}_{\G}:=X_{\G}\cup\lbrace \infty\rbrace$ denotes the one-point compactfication of $X_{\G}$. For any $k\geq 1$, let $\lambda^{k}:=\Lambda^{\otimes k}\otimes O_{\overline{X}_{\G}}((k-1)\infty)$. For any $k\geq 1$, the hyperbolic metric induces a point-wise metric on $H^{0}(\overline{X}_{\G},\lambda^{k})$, which we denote by $|\cdot|_{\mathrm{hyp}}$. For any $k\geq 1$, let $\bkx$ denote the Bergman kernel associated to the complex vector space of $H^{0}(\overline{X}_{\G},\lambda^{k})$. Then, for $k\gg1$, the first main result of the article, is the following estimate
\begin{align*}
\sup_{z\in \overline{X}_{\G}}\big| \bkx(z,z)\big|_{\mathrm{hyp}}=O_{X_{\G}}\big(k^{n+1\slash 2} \big),
\end{align*}
%%%%%%%%%%%%%%%%%%%%%%%%%%%%%%%%%%%%%%%%%%%%%%%%%%%%%%%%%%%%%%%%%%%%%%%%%%%%%%%%%%%%%%%%%%%
where the implied constant depends only on $X_{\G}$. 
%%%%%%%%%%%%%%%%%%%%%%%%%%%%%%%%%%%%%%%%%%%%%%%%%%%%%%%%%%%%%%%%%%%%%%%%%%%%%%%%%%%%%%%%%%%
 
\vspace{0.1cm}\noindent
For any $k\geq 1$, and $z\in X_{\G}$, let 
\begin{align*}
\mu_{\mathrm{Ber},k}(z):=-\frac{i}{2\pi}\partial_{z}\partial_{\overline{z}}\big| \bkx(z,z)\big|_{\mathrm{hyp}}.
\end{align*}
%%%%%%%%%%%%%%%%%%%%%%%%%%%%%%%%%%%%%%%%%%%%%%%%%%%%%%%%%%%%%%%%%%%%%%%%%%%%%%%%%%%%%%%%%%%
denote the Bergman metric associated to the line bundle $\lambda^{ k}$, and let $\mu_{\mathrm{Ber},k}^{\mathrm{vol}}$ denote the associated volume form.  Then, for $k\gg1$, the second main result of the article is the following estimate
\begin{align*}
\sup_{z\in \overline{X}_{\G}}\bigg|\frac{\mu_{\mathrm{Ber},k}^{\mathrm{vol}}(z)}{\mu_{\mathrm{hyp}}^{\mathrm{vol}}(z)}\bigg|=O_{X_{\G}}\big(k^{2(n-1)(n+2)+n+3} \big),
\end{align*}
%%%%%%%%%%%%%%%%%%%%%%%%%%%%%%%%%%%%%%%%%%%%%%%%%%%%%%%%%%%%%%%%%%%%%%%%%%%%%%%%%%%%%%%%%%%
where the implied constant depends only on $X_{\G}$.
\end{abstract}
%%%%%%%%%%%%%%%%%%%%%%%%%%%%%%%%%%%%%%%%%%%%%%%%%%%%%%%%%%%%%%%%%%%%%%%%%%%%%%%%%%%%%%%%%%%
%%%%%%%%%%%%%%%%%%%%%%%%%%%%%%%%%%%%%%%%%%%%%%%%%%%%%%%%%%%%%%%%%%%%%%%%%%%%%%%%%%%%%%%%%%%
\maketitle
%%%%%%%%%%%%%%%%%%%%%%%%%%%%%%%%%%%%%%%%%%%%%%%%%%%%%%%%%%%%%%%%%%%%%%%%%%%%%%%%%%%%%%%%%%%
\section{Introduction}\label{sec-1}
%%%%%%%%%%%%%%%%%%%%%%%%%%%%%%%%%%%%%%%%%%%%%%%%%%%%%%%%%%%%%%%%%%%%%%%%%%%%%%%%%%%%%%%%%%%
In section \ref{subsec1.1}, we discuss the history and motivation behind the results proved, and in section \ref{subsec1.2}, we state the main results proved in this article. 
%%%%%%%%%%%%%%%%%%%%%%%%%%%%%%%%%%%%%%%%%%%%%%%%%%%%%%%%%%%%%%%%%%%%%%%%%%%%%%%%%%%%%%%%%%%

\vspace{0.2cm}
\subsection{History}\label{subsec1.1}
%%%%%%%%%%%%%%%%%%%%%%%%%%%%%%%%%%%%%%%%%%%%%%%%%%%%%%%%%%%%%%%%%%%%%%%%%%%%%%%%%%%%%%%%%%%
Estimates of Bergman kernels associated to high tensor powers of holomorphic line bundles defined over complex manifolds, is a well studied problem in Complex Geometry, by the likes of Tian, Zeldtich, Demailly, Catalin, Ma, and Marinescu, et al. Especially, estimates of Bergman kernels associated to high tensor powers of the cotangent bundle, in the context of arithmetic manifolds, are of great importance in deriving individual estimates of Hecke eigen automorphic forms. Deriving such estimates is an area of immense interest from the perspective of analytic number theory, as well as arithmetic intersection theory, which we now describe.  
%%%%%%%%%%%%%%%%%%%%%%%%%%%%%%%%%%%%%%%%%%%%%%%%%%%%%%%%%%%%%%%%%%%%%%%%%%%%%%%%%%%%%%%%%%%

\vspace{0.1cm}
In \cite{auvray} and \cite{jk}, following different methods, the authors have derived optimal estimates of the Bergman kernel associated to the cotangent bundle, defined over compact and nonconpact hyperbolic Riemann surfaces of finite hyperbolic volume. Similar estimates were derived in \cite{anil1} and \cite{anil2}. In \cite{anil3}, the authors derived optimal estimates of the Bergman kernel associated to the cotangent bundle, defined over compact and nonconpact Hilbert modular varieties of finite hyperbolic volume. In \cite{das}, the authors have extended the estimates from \cite{jk}, to the setting of Siegel three-folds. In \cite{mandal}, the authors has extended the estimates from \cite{jk} to Siegel modular varieties of arbitrary dimension. 
%%%%%%%%%%%%%%%%%%%%%%%%%%%%%%%%%%%%%%%%%%%%%%%%%%%%%%%%%%%%%%%%%%%%%%%%%%%%%%%%%%%%%%%%%%%

\vspace{0.1cm}
In \cite{anil4}, the authors have extended the estimates from \cite{jk} to Picard surfaces. From evidence emanating from the estimates derived in \cite{jk}, \cite{anil3}, \cite{das}, and \cite{mandal}, a folk-lore conjecture was expected to hold. However, estimates from \cite{anil4} are stronger than the estimate predicted by the folk-lore conjecture.  
%%%%%%%%%%%%%%%%%%%%%%%%%%%%%%%%%%%%%%%%%%%%%%%%%%%%%%%%%%%%%%%%%%%%%%%%%%%%%%%%%%%%%%%%%%%

\vspace{0.1cm}
The first main result of the article is an extension of the estimates derived in \cite{anil4} to the setting of Picard varieties of arbitrary dimension. Our estimate is stronger than the  conjectured estimate in \cite{anil4}, and is optimally derived.  
%%%%%%%%%%%%%%%%%%%%%%%%%%%%%%%%%%%%%%%%%%%%%%%%%%%%%%%%%%%%%%%%%%%%%%%%%%%%%%%%%%%%%%%%%%%

\vspace{0.1cm}
Using a different set of techniques, a similar estimate as our first main result was derived by Zhou in \cite{zhou}. The approach in \cite{zhou} is an application of the method of peak sections by Tian.
%%%%%%%%%%%%%%%%%%%%%%%%%%%%%%%%%%%%%%%%%%%%%%%%%%%%%%%%%%%%%%%%%%%%%%%%%%%%%%%%%%%%%%%%%%%

\vspace{0.1cm}
Another problem of interest in complex geometry, is deriving estimates of Bergman metric, associated to tensor powers of a given holomorphic line bundle. In \cite{anil5}, Biswas and the first named author have derived estimates of the Bergman metric associated to high tensor-powers of the cotangent bundle, defined over a compact hyperbolic Riemann surface. In \cite{anil6}, the first named author and Mukherjee have extended estimates from \cite{anil5} to the setting of noncompact hyperbolic Riemann surfaces. In \cite{anil7}, the first named author, Roy, and Sadhukhan have extended the estimates from \cite{anil5} to the setting of Picard surfaces. 
%%%%%%%%%%%%%%%%%%%%%%%%%%%%%%%%%%%%%%%%%%%%%%%%%%%%%%%%%%%%%%%%%%%%%%%%%%%%%%%%%%%%%%%%%%%

\vspace{0.1cm}
Estimates of Bergman metric in the setting of noncompact manifolds are difficult to obtain, as the Bergman metric and the ambient metric on the manifold, both admit singularities. Improving and improvising arguments from \cite{anil7}, we derive  estimates of the Bergman metric associated to high tensor-powers of the cotangent bundle derived in \cite{anil5} and \cite{anil7}, to the setting of noncompact Picard varieties, of finite hyperbolic volume, which is the second main result of this article. 
%%%%%%%%%%%%%%%%%%%%%%%%%%%%%%%%%%%%%%%%%%%%%%%%%%%%%%%%%%%%%%%%%%%%%%%%%%%%%%%%%%%%%%%%%%%

\vspace{0.2cm}
\subsection{Statement of results}\label{subsec1.2}
%%%%%%%%%%%%%%%%%%%%%%%%%%%%%%%%%%%%%%%%%%%%%%%%%%%%%%%%%%%%%%%%%%%%%%%%%%%%%%%%%%%%%%%%%%%
We now describe the two main results of this article. 
%%%%%%%%%%%%%%%%%%%%%%%%%%%%%%%%%%%%%%%%%%%%%%%%%%%%%%%%%%%%%%%%%%%%%%%%%%%%%%%%%%%%%%%%%%%

\vspace{0.1cm}
For a fixed $n\geq 2$, let $\G\subset\mathrm{SU}((n,1),\calo_{K})$ be a torsion-free, finite-index subgroup, with only one cusp at $\infty$, where $K$ is a totally imaginary number field of degree $2$. Let $\Bn$ denote the unit complex ball endowed with the hyperbolic metric $\hyp$, and we denote the volume form associated to the hyperbolic metric by $\hypbnvol$. Let $X_{\G}:=\G\backslash \Bn$ denote the quotient space, which is a complex hyperbolic manifold of dimension $n$, and is of finite hyperbolic volume. Furthermore, let 
$\overline{X}_{\G}:=X_{\G}\cup \lbrace \infty\rbrace$ denote the one-point compactification of $X_{\G}$. 
%%%%%%%%%%%%%%%%%%%%%%%%%%%%%%%%%%%%%%%%%%%%%%%%%%%%%%%%%%%%%%%%%%%%%%%%%%%%%%%%%%%%%%%%%%%

\vspace{0.1cm}
Let $\Lambda:=\Omega^{n}_{\overline{X}_{\G}}$ denote the line bundle of holomorphic $n$-forms, defined over $\overline{X}_{\G}$. For any $k\geq 1$, let $\lambda^{k}:=\Lambda^{\otimes k}\otimes O_{\overline{X}_{\G}}((k-1)\infty)$, and let $\bkx$ denote the Bergman kernel associated to  $H^{0}(\overline{X}_{\G}, \lambda^{k})$, the complex vector space of holomorphic global sections of the line bundle $\lambda^{ k}$. The hyperbolic metric induces a point-wise metric on $H^{0}(\overline{X}_{\G}, \lambda^{ k})$, which we denote by $|\cdot|_{\mathrm{hyp}}$. 
%%%%%%%%%%%%%%%%%%%%%%%%%%%%%%%%%%%%%%%%%%%%%%%%%%%%%%%%%%%%%%%%%%%%%%%%%%%%%%%%%%%%%%%%%%%

\vspace{0.1cm}
The first main result of the article is the following theorem, which is proved as Theorem \ref{thm4} in section \ref{subsec-3.1}. 
%%%%%%%%%%%%%%%%%%%%%%%%%%%%%%%%%%%%%%%%%%%%%%%%%%%%%%%%%%%%%%%%%%%%%%%%%%%%%%%%%%%%%%%%%%%

\vspace{0.1cm}
\begin{mainthm}\label{mainthm1}
With notation as above, for $k\gg 1$, we have the following estimate
\begin{align}\label{est:mainthm1}
\sup_{z\in \overline{X}_{\G}}\big| \bkx(z,z)\big|_{\mathrm{hyp}}=O_{X_{\G}}\big(k^{n+1\slash 2}\big),
\end{align}
%%%%%%%%%%%%%%%%%%%%%%%%%%%%%%%%%%%%%%%%%%%%%%%%%%%%%%%%%%%%%%%%%%%%%%%%%%%%%%%%%%%%%%%%%%%
where the implied constant depends only on $X_{\G}$.
\end{mainthm}
%%%%%%%%%%%%%%%%%%%%%%%%%%%%%%%%%%%%%%%%%%%%%%%%%%%%%%%%%%%%%%%%%%%%%%%%%%%%%%%%%%%%%%%%%%%

\vspace{0.1cm}
For a fixed $n\geq 2$, and $k\geq 1$, the Bergman metric associated to the line bundle $\Lambda^{\otimes k}$, is given by the following formula
\begin{align*}
\mu_{\mathrm{Ber},k}(z):=-\frac{i}{2\pi}\partial_{z}\partial_{\overline{z}}\log\big|\bkx(z,z)\big|_{\mathrm{hyp}}.
\end{align*}
%%%%%%%%%%%%%%%%%%%%%%%%%%%%%%%%%%%%%%%%%%%%%%%%%%%%%%%%%%%%%%%%%%%%%%%%%%%%%%%%%%%%%%%%%%%
Let $\mu_{\mathrm{Ber},k}^{\mathrm{vol}}$ denote the volume form associated the metric $\mu_{\mathrm{Ber},k}$. 
%%%%%%%%%%%%%%%%%%%%%%%%%%%%%%%%%%%%%%%%%%%%%%%%%%%%%%%%%%%%%%%%%%%%%%%%%%%%%%%%%%%%%%%%%%%

\vspace{0.1cm}
The second main result of the article is the following theorem, which is proved as Theorem \ref{thm8} in section \ref{subsec-3.2}. 
%%%%%%%%%%%%%%%%%%%%%%%%%%%%%%%%%%%%%%%%%%%%%%%%%%%%%%%%%%%%%%%%%%%%%%%%%%%%%%%%%%%%%%%%%%%

\vspace{0.1cm}
\begin{mainthm}\label{mainthm2}
With notation as above, for $k\gg1$, we have the following estimate
\begin{align}\label{est:mainthm2}
\sup_{z\in \overline{X}_{\G}}\bigg|\frac{\mu_{\mathrm{Ber},k}^{\mathrm{vol}}(z)}{\hypbnvol(z)}\bigg|=O_{X_{\G}}\big(k^{2(n-1)(n+2)+n+3} \big),
\end{align}
%%%%%%%%%%%%%%%%%%%%%%%%%%%%%%%%%%%%%%%%%%%%%%%%%%%%%%%%%%%%%%%%%%%%%%%%%%%%%%%%%%%%%%%%%%%
where the implied constant depends only on $X_{\G}$.
\end{mainthm}
%%%%%%%%%%%%%%%%%%%%%%%%%%%%%%%%%%%%%%%%%%%%%%%%%%%%%%%%%%%%%%%%%%%%%%%%%%%%%%%%%%%%%%%%%%%

\vspace{0.1cm}
\begin{rem}
The restriction on the number of cusps, can be relaxed, via scaling matrices. For the brevity of notation, and of exposition, we restrict ourselves to setting of a single cusp at $\infty$.   
%%%%%%%%%%%%%%%%%%%%%%%%%%%%%%%%%%%%%%%%%%%%%%%%%%%%%%%%%%%%%%%%%%%%%%%%%%%%%%%%%%%%%%%%%%%

\vspace{0.1cm}
As per the folk-lore conjecture, the following estimate is expected to hold true
\begin{align}\label{folk-conj}
\sup_{z\in \overline{X}_{\G}}\big| \bkx(z,z)\big|_{\mathrm{hyp}}=O_{X_{\G}}\big(k^{3n\slash 2}\big), 
\end{align}
%%%%%%%%%%%%%%%%%%%%%%%%%%%%%%%%%%%%%%%%%%%%%%%%%%%%%%%%%%%%%%%%%%%%%%%%%%%%%%%%%%%%%%%%%%%
where the implied constant depends only on $X_{\G}$. 
%%%%%%%%%%%%%%%%%%%%%%%%%%%%%%%%%%%%%%%%%%%%%%%%%%%%%%%%%%%%%%%%%%%%%%%%%%%%%%%%%%%%%%%%%%%

\vspace{0.1cm}
Estimate  \eqref{est:mainthm1} is much stronger than estimate \eqref{folk-conj}, and is stronger than  conjecture made in \cite{anil4}.
%%%%%%%%%%%%%%%%%%%%%%%%%%%%%%%%%%%%%%%%%%%%%%%%%%%%%%%%%%%%%%%%%%%%%%%%%%%%%%%%%%%%%%%%%%%
\end{rem}
%%%%%%%%%%%%%%%%%%%%%%%%%%%%%%%%%%%%%%%%%%%%%%%%%%%%%%%%%%%%%%%%%%%%%%%%%%%%%%%%%%%%%%%%%%%

\vspace{0.1cm}
\begin{rem}
Our estimate of the Bergman metric \eqref{est:mainthm2} is far from being optimal. The power of $k$ in our estimate \eqref{est:mainthm2} is quadratic in $n$, and we believe it should linear in $n$. For $n\geq 1$, we conjecture the following estimate
\begin{align*}
\sup_{z\in \overline{X}_{\G}}\bigg|\frac{\mu_{\mathrm{Ber},k}^{\mathrm{vol}}(z)}{\hypbnvol(z)}\bigg|=O_{X_{\G}}\big(k^{4n}\big),
\end{align*}
%%%%%%%%%%%%%%%%%%%%%%%%%%%%%%%%%%%%%%%%%%%%%%%%%%%%%%%%%%%%%%%%%%%%%%%%%%%%%%%%%%%%%%%%%%%
where the implied constant depends only on $X_{\G}$. 
%%%%%%%%%%%%%%%%%%%%%%%%%%%%%%%%%%%%%%%%%%%%%%%%%%%%%%%%%%%%%%%%%%%%%%%%%%%%%%%%%%%%%%%%%%%

\vspace{0.1cm}
As per evidence from the one-dimensional setting from \cite{amw},  our estimate \eqref{est:mainthm2} is optimal for the case $n=1$. It is to be mentioned that the lower bound proved in Proposition \ref{prop9}, corrects and completes arguments proved in \cite{anil6}, for the case $n=1$. Hence, for $n=1$, we have the following estimate
\begin{align*}
\sup_{z\in \overline{X}_{\G}}\bigg|\frac{\mu_{\mathrm{Ber},k}^{\mathrm{vol}}(z)}{\hypbnvol(z)}\bigg|=O_{X_{\G}}\big(k^{4}\big),
\end{align*}
%%%%%%%%%%%%%%%%%%%%%%%%%%%%%%%%%%%%%%%%%%%%%%%%%%%%%%%%%%%%%%%%%%%%%%%%%%%%%%%%%%%%%%%%%%%
\end{rem}
%%%%%%%%%%%%%%%%%%%%%%%%%%%%%%%%%%%%%%%%%%%%%%%%%%%%%%%%%%%%%%%%%%%%%%%%%%%%%%%%%%%%%%%%%%%
%%%%%%%%%%%%%%%%%%%%%%%%%%%%%%%%%%%%%%%%%%%%%%%%%%%%%%%%%%%%%%%%%%%%%%%%%%%%%%%%%%%%%%%%%%%

\vspace{0.2cm}
\section{Background material}\label{sec-2}
%%%%%%%%%%%%%%%%%%%%%%%%%%%%%%%%%%%%%%%%%%%%%%%%%%%%%%%%%%%%%%%%%%%%%%%%%%%%%%%%%%%%%%%%%%%
In this section, we set up the notation and state the relevant results, used to prove Main Theorem \ref{mainthm1} and Main Theorem \ref{mainthm2}.  
%%%%%%%%%%%%%%%%%%%%%%%%%%%%%%%%%%%%%%%%%%%%%%%%%%%%%%%%%%%%%%%%%%%%%%%%%%%%%%%%%%%%%%%%%%%

\vspace{0.2cm}
\subsection{The ball model}\label{subsec-2.1}
%%%%%%%%%%%%%%%%%%%%%%%%%%%%%%%%%%%%%%%%%%%%%%%%%%%%%%%%%%%%%%%%%%%%%%%%%%%%%%%%%%%%%%%%%%%
For $n\geq 2$, let  $H$  be a Hermitian matrix of signature $(n,1)$, and let $ \big(\mathbb{C}^{n+1},H\big)$ denote the Hermitian inner-product space $\C^{n+1}$ equipped with the Hermitian inner-product $\langle\cdot, \cdot\rangle_{H}$, induced by $H$. For $z,w \in \mathbb{C}^{n+1}$, the Hermitian inner-product $\langle\cdot, \cdot\rangle_{H}$ is defined as 
\begin{align*}
\langle z, w \rangle_H = w^{*}H z,
\end{align*}
%%%%%%%%%%%%%%%%%%%%%%%%%%%%%%%%%%%%%%%%%%%%%%%%%%%%%%%%%%%%%%%%%%%%%%%%%%%%%%%%%%%%%%%%%%%
where $w^{*}$ denotes the complex transpose of $w$.
%%%%%%%%%%%%%%%%%%%%%%%%%%%%%%%%%%%%%%%%%%%%%%%%%%%%%%%%%%%%%%%%%%%%%%%%%%%%%%%%%%%%%%%%%%%

\vspace{0.1cm}
Set
\begin{align*}
V_{-}(H):= \{z\in\mathbb{C}^{n+1}\big|\,\langle z, z \rangle_H < 0 \}\subset \C^{n+1};\,\,
 V_0(H):= \{z\in\mathbb{C}^{n+1}\big|\,\langle z, z \rangle_H = 0 \}\subset \C^{n+1}.
\end{align*}
%%%%%%%%%%%%%%%%%%%%%%%%%%%%%%%%%%%%%%%%%%%%%%%%%%%%%%%%%%%%%%%%%%%%%%%%%%%%%%%%%%%%%%%%%%%

\vspace{0.1cm}
Set 
\begin{align*}
\cala^{n+1}:=\big\lbrace \big(z_{1},\ldots,z_{n+1}\big)^{t}\in \C^{n+1}\big|\,z_{n+1}\not=0\rbrace\subset \C^{n+1}.
\end{align*}
%%%%%%%%%%%%%%%%%%%%%%%%%%%%%%%%%%%%%%%%%%%%%%%%%%%%%%%%%%%%%%%%%%%%%%%%%%%%%%%%%%%%%%%%%%%

\vspace{0.1cm}
The subspace $\cala^{n+1}$ inherits the structure of a Hermitian inner-product space, which we denote by $\big(\cala^{n+1},H\big)$. We have the following map
\begin{align*}
&\phi:\cala_{n+1} \longrightarrow \C^{n}\\&
(z_1,\ldots,z_{n+1})^{t}\mapsto \bigg( \frac{z_1}{z_{n+1}},  \cdots ,\frac{z_n}{z_{n+1}}\bigg)^{t}.  
\end{align*}                    
%%%%%%%%%%%%%%%%%%%%%%%%%%%%%%%%%%%%%%%%%%%%%%%%%%%%%%%%%%%%%%%%%%%%%%%%%%%%%%%%%%%%%%%%%%%

\vspace{0.1cm}
The complex hyperbolic space  $\mathcal{H}^{n}\left(\mathbb{C}\right)$  and its boundary $\partial \mathcal{H}^{n}\left(\mathbb{C}\right)$ are defined as 
\begin{align*}
\mathcal{H}^{n}\left(\mathbb{C}\right):=\phi\big( V_{-}(H)\big)\subset \C^{n};\,\,\,\partial \mathcal{H}^{n}\left(\mathbb{C}\right):=\phi\big( V_{0}(H)\big)\subset\C^{n}.
\end{align*}
%%%%%%%%%%%%%%%%%%%%%%%%%%%%%%%%%%%%%%%%%%%%%%%%%%%%%%%%%%%%%%%%%%%%%%%%%%%%%%%%%%%%%%%%%%%

\vspace{0.1cm}
The choice of Hermitian matrix gives different models of hyperbolic $n$-space. For
\begin{align}\label{defn:H}
H = \begin{pmatrix}
\mathrm{Id}_{n } & 0 \\
0  & -1
\end{pmatrix},
\end{align}
%%%%%%%%%%%%%%%%%%%%%%%%%%%%%%%%%%%%%%%%%%%%%%%%%%%%%%%%%%%%%%%%%%%%%%%%%%%%%%%%%%%%%%%%%%%
we obtain the unit ball model of the complex hyperbolic space, where $\mathrm{Id}_{n}$ denotes the identity matrix of order $n$, which we now describe.
%%%%%%%%%%%%%%%%%%%%%%%%%%%%%%%%%%%%%%%%%%%%%%%%%%%%%%%%%%%%%%%%%%%%%%%%%%%%%%%%%%%%%%%%%%%

\vspace{0.1cm}
For any $z= (z_1,\ldots ,z_n)^{t} \in \mathcal{H}^n(\C)$, consider the lift $\tilde{z}:=\big(z_1,\ldots ,z_n,1\big)^{t}\in\cala^{n+1}$. Observe that
\begin{align*}
\langle \tilde{z},\tilde{z}\rangle_{H} < 0
\iff  |z_1|^2 +\cdots +|z_n|^2 <1.
\end{align*}
%%%%%%%%%%%%%%%%%%%%%%%%%%%%%%%%%%%%%%%%%%%%%%%%%%%%%%%%%%%%%%%%%%%%%%%%%%%%%%%%%%%%%%%%%%%

Thus $\mathcal{H}^n(\C)$ can be identified with the open unit ball 
\begin{align*}
\mathbb{B}^n := \big\lbrace z:=(z_{1},\ldots,z_{n})^{t}\in \C^{n}\big|\, |z|^{2}:=|z_1|^2+ \cdots + |z_n|^2 < 1 \big\rbrace.
\end{align*}
%%%%%%%%%%%%%%%%%%%%%%%%%%%%%%%%%%%%%%%%%%%%%%%%%%%%%%%%%%%%%%%%%%%%%%%%%%%%%%%%%%%%%%%%%%%
This model which is also known as the ball model, is the natural generalization of the Poincare disk model of the upper half-plane $\mathbb{H}$.
%%%%%%%%%%%%%%%%%%%%%%%%%%%%%%%%%%%%%%%%%%%%%%%%%%%%%%%%%%%%%%%%%%%%%%%%%%%%%%%%%%%%%%%%%%%
 
\vspace{0.1cm} 
The space is endowed with its usual K\"ahler-Bergman metric $\mu_{\mathrm{hyp}}$,  which has constant negative curvature $-1$. For any given $z=(z_{1},\ldots,z_{n})\in \mathbb{B}^n$, the hyperbolic metric $\mu_{\mathrm{hyp}}(z)$  is given by the following formula
\begin{align}\label{defnhypmetric}
\mu_{\mathrm{hyp}}(z):=-\frac{i}{2\pi}\partial\overline{\partial}\log(1-|z|^{2}),
\end{align}
%%%%%%%%%%%%%%%%%%%%%%%%%%%%%%%%%%%%%%%%%%%%%%%%%%%%%%%%%%%%%%%%%%%%%%%%%%%%%%%%%%%%%%%%%%%
and let $\mu_{\mathrm{hyp}}^{\mathrm{vol}}(z)$ denote the associated volume form, which is given by the following formula
\begin{align}\label{defnhypmetric-vol}
\mu_{\mathrm{hyp}}^{\mathrm{vol}}(z):=\bigg(\frac{i}{2\pi}\bigg)^{n}\frac{dz_1\wedge d\overline{z}_{1}\wedge\cdots\wedge dz_n\wedge d\overline{z}_{n} }{\big( 1-|z|^{2}\big)^{n+1}}
\end{align}
%%%%%%%%%%%%%%%%%%%%%%%%%%%%%%%%%%%%%%%%%%%%%%%%%%%%%%%%%%%%%%%%%%%%%%%%%%%%%%%%%%%%%%%%%%%

\vspace{0.1cm}
From standard results in hyperbolic geometry, for any $z=(z_1,\ldots,z_n)^{t},\,w=(w_1,\ldots,w_n)^{t}\in \mathbb{B}^n$ with respective lifts $\tilde{z}=(z_1,\ldots,z_n,1)^{t},\,\tilde{w}=(w_1,\ldots,w_n,1)^{t}\in \cala^{n+1}$, the hyperbolic distance in the ball model is given by the following relation
\begin{align}\label{coshreln}
&\cosh^{2}\big(\dhyp(z,w)\slash 2\big)=\frac{\langle \tilde{z},\tilde{w}\rangle_H\; \langle \tilde{w},\tilde{z} \rangle_H }{\langle \tilde{z},\tilde{z}\rangle_H\; \langle \tilde{w},\tilde{w}\rangle_H},\,\,\,\;\mathrm{where}\,\,\langle \tilde{z},\tilde{w} \rangle _H:=\tilde{w}^{\ast}H\tilde{z}=\sum_{j=1}^{n}z_j\overline{w}_{j}-1,
\end{align}
%%%%%%%%%%%%%%%%%%%%%%%%%%%%%%%%%%%%%%%%%%%%%%%%%%%%%%%%%%%%%%%%%%%%%%%%%%%%%%%%%%%%%%%%%%%
and $H$ is as described in equation \eqref{defn:H}.
%%%%%%%%%%%%%%%%%%%%%%%%%%%%%%%%%%%%%%%%%%%%%%%%%%%%%%%%%%%%%%%%%%%%%%%%%%%%%%%%%%%%%%%%%%%
%%%%%%%%%%%%%%%%%%%%%%%%%%%%%%%%%%%%%%%%%%%%%%%%%%%%%%%%%%%%%%%%%%%%%%%%%%%%%%%%%%%%%%%%%%%

\vspace{0.2cm} 
\subsection{Alternate models of $\calh^n (\C)$}\label{subsec-2.2}
%%%%%%%%%%%%%%%%%%%%%%%%%%%%%%%%%%%%%%%%%%%%%%%%%%%%%%%%%%%%%%%%%%%%%%%%%%%%%%%%%%%%%%%%%%%
Different choices of $H$ give us different models for $\calh^n (\C)$. When $H$ is taken as the matrix $H_1 = \begin{pmatrix} I_n & \\ & -1 \end{pmatrix}$, we get the ball model considered earlier, and we will refer to this as model (1).
%%%%%%%%%%%%%%%%%%%%%%%%%%%%%%%%%%%%%%%%%%%%%%%%%%%%%%%%%%%%%%%%%%%%%%%%%%%%%%%%%%%%%%%%%%%

\vspace{0.1cm} 
When we take $H$ to be the matrix
\begin{align*}
H_2 = \begin{pmatrix} & & i \\ & I_{n-1} & \\  -i & & \end{pmatrix},
\end{align*}
%%%%%%%%%%%%%%%%%%%%%%%%%%%%%%%%%%%%%%%%%%%%%%%%%%%%%%%%%%%%%%%%%%%%%%%%%%%%%%%%%%%%%%%%%%%
we get the upper unbounded hyperquadratic model of the complex hyperbolic space. This space is given by 
\begin{align*}
\D^n = \bigg\lbrace z:=(z_{1},\ldots,z_{n})^{t}\in \C^{n}\big|\, \mathrm{Im}(z_n) > \frac{1}{2} \sum_{j=1}^{n-1} |z_j|^2 \bigg\rbrace,
\end{align*}
%%%%%%%%%%%%%%%%%%%%%%%%%%%%%%%%%%%%%%%%%%%%%%%%%%%%%%%%%%%%%%%%%%%%%%%%%%%%%%%%%%%%%%%%%%%
which we denote as model (2).
%%%%%%%%%%%%%%%%%%%%%%%%%%%%%%%%%%%%%%%%%%%%%%%%%%%%%%%%%%%%%%%%%%%%%%%%%%%%%%%%%%%%%%%%%%%

\vspace{0.1cm} 
Finally, by taking $H$ to be the matrix
\begin{align}\label{defn:H_3}
H_3 = \begin{pmatrix} & & 1 \\ & I_{n-1} & \\  1 & & \end{pmatrix},
\end{align}
%%%%%%%%%%%%%%%%%%%%%%%%%%%%%%%%%%%%%%%%%%%%%%%%%%%%%%%%%%%%%%%%%%%%%%%%%%%%%%%%%%%%%%%%%%%
we get the left unbounded hyperquadratic model of the complex hyperbolic space.  This space is given by 
\begin{align*}
\E^n = \bigg\lbrace z:=(z_{1},\ldots,z_{n})^{t}\in \C^{n}\big|\, \mathrm{Re}(z_1) < -\frac{1}{2} \sum_{j=2}^{n} |z_j|^2 \bigg\rbrace,
\end{align*}
%%%%%%%%%%%%%%%%%%%%%%%%%%%%%%%%%%%%%%%%%%%%%%%%%%%%%%%%%%%%%%%%%%%%%%%%%%%%%%%%%%%%%%%%%%%
which we denote as model (3). We wish to be able to go back and forth between these models. We define the following three Cayley transforms, which are all isometries.
%%%%%%%%%%%%%%%%%%%%%%%%%%%%%%%%%%%%%%%%%%%%%%%%%%%%%%%%%%%%%%%%%%%%%%%%%%%%%%%%%%%%%%%%%%%

\vspace{0.1cm} 
We have the following isometry between the models (1) and (2). Consider the map
\begin{align*}
\gamma_{21} : (\D^n, H_2)& \longrightarrow (\mathbb{B}^n, H_1)\\
(z_1, \dots, z_n)& \mapsto \left( \frac{z_1}{z_n + i}, \dots, \frac{z_{n-1}}{z_n + i}, \frac{z_n - i}{z_n + i} \right).
\end{align*}
%%%%%%%%%%%%%%%%%%%%%%%%%%%%%%%%%%%%%%%%%%%%%%%%%%%%%%%%%%%%%%%%%%%%%%%%%%%%%%%%%%%%%%%%%%%
The definition of $z \in \D^n$ can be rephrased as $-2 \mathrm{Im} (z_n) < - \sum_{j=1}^{n-1} |z_j|^2$. This ensures that $\gamma_{21} (z) \in \mathbb{B}^n$.
%%%%%%%%%%%%%%%%%%%%%%%%%%%%%%%%%%%%%%%%%%%%%%%%%%%%%%%%%%%%%%%%%%%%%%%%%%%%%%%%%%%%%%%%%%%

\vspace{0.1cm} 
Similarly, we define the following isometry between the models (1) and (3). Consider the map
\begin{align}\label{gamma31}
\gamma_{31} : (\E^n, H_3) &\longrightarrow (\mathbb{B}^n, H_1)\notag\\[0.1cm]
(z_1, \dots, z_n) &\mapsto \left( \frac{z_1 + 1}{z_1 - 1},  \frac{\sqrt{2}z_2}{z_1 - 1}, \dots, \frac{\sqrt{2}z_{n}}{z_1 - 1} \right).
\end{align}
%%%%%%%%%%%%%%%%%%%%%%%%%%%%%%%%%%%%%%%%%%%%%%%%%%%%%%%%%%%%%%%%%%%%%%%%%%%%%%%%%%%%%%%%%%%
It is again a straightforward check that for $z \in \E^n$, we have $\gamma_{31} (z) \in \mathbb{B}^n$, and the inverse map is given by
\begin{align}\label{gamma13}
\gamma_{13}: (\mathbb{B}^n, H_1) &\longrightarrow (\mathbb{E}^n, H_3)\notag\\[0.1cm]
(z_1, \dots, z_n) &\mapsto \left( \frac{z_1 + 1}{z_1 - 1},  \frac{\sqrt{2}z_2}{z_1 - 1}, \dots, \frac{\sqrt{2}z_{n}}{z_1 - 1} \right).
\end{align}
%%%%%%%%%%%%%%%%%%%%%%%%%%%%%%%%%%%%%%%%%%%%%%%%%%%%%%%%%%%%%%%%%%%%%%%%%%%%%%%%%%%%%%%%%%%

\vspace{0.1cm} 
Finally we define $\gamma_{23} = \gamma_{13} \circ \gamma_{21}: (\D^n, H_2) \to (\E^n, H_3)$, which is an isometry between models (2) and (3). 
%%%%%%%%%%%%%%%%%%%%%%%%%%%%%%%%%%%%%%%%%%%%%%%%%%%%%%%%%%%%%%%%%%%%%%%%%%%%%%%%%%%%%%%%%%%

\vspace{0.1cm}
From standard results in hyperbolic geometry, for any $z=(z_1,\ldots,z_n)^{t},\,w=(w_1,\ldots,w_n)^{t}\in \E^n$ with respective lifts $\tilde{z}=(z_1,\ldots,z_n,1)^{t},\,\tilde{w}=(w_1,\ldots,w_n,1)^{t}\in \cala^{n+1}$, the hyperbolic distance in model (3), is given by the following relation
\begin{align}\label{coshreln-1}
&\cosh^{2}\big(\dhyp(z,w)\slash 2\big)=\frac{\langle \tilde{z},\tilde{w}\rangle_{H_3}\; \langle \tilde{w},\tilde{z} \rangle_{H_3} }{\langle \tilde{z},\tilde{z}\rangle_{H_3}\; \langle \tilde{w},\tilde{w}\rangle_{H_{3}}},\notag\\[0.12cm]&\mathrm{where}\,\,\langle \tilde{z},\tilde{w} \rangle_{H_3}:=\tilde{w}^{\ast}H_{3}\tilde{z}=z_1+\overline{w_{1}}+\sum_{j=2}^{n}z_j\overline{w}_{j},
\end{align}
%%%%%%%%%%%%%%%%%%%%%%%%%%%%%%%%%%%%%%%%%%%%%%%%%%%%%%%%%%%%%%%%%%%%%%%%%%%%%%%%%%%%%%%%%%%
and $H_3$ is as described in equation \eqref{defn:H_3}.
%%%%%%%%%%%%%%%%%%%%%%%%%%%%%%%%%%%%%%%%%%%%%%%%%%%%%%%%%%%%%%%%%%%%%%%%%%%%%%%%%%%%%%%%%%%
%%%%%%%%%%%%%%%%%%%%%%%%%%%%%%%%%%%%%%%%%%%%%%%%%%%%%%%%%%%%%%%%%%%%%%%%%%%%%%%%%%%%%%%%%%%

\vspace{0.2cm} 
\subsection{Picard variety}\label{subsec-2.3}
%%%%%%%%%%%%%%%%%%%%%%%%%%%%%%%%%%%%%%%%%%%%%%%%%%%%%%%%%%%%%%%%%%%%%%%%%%%%%%%%%%%%%%%%%%%
The unitary group $\mathrm{SU}((n,1),\C)$ is given by 
\begin{align*}
\mathrm{SU}((n,1),\C):= \{ A \in \mathrm{SL}(n+1,\C) \; | \, A^*H A =  H\},\, \,\mathrm{where} \,\,
H = \begin{pmatrix}
\mathrm{Id}_{n } & 0 \\
0  & -1
\end{pmatrix},
\end{align*}
%%%%%%%%%%%%%%%%%%%%%%%%%%%%%%%%%%%%%%%%%%%%%%%%%%%%%%%%%%%%%%%%%%%%%%%%%%%%%%%%%%%%%%%%%%%
and $\mathrm{Id}_{n}$ denotes the identity matrix of order $n$.
%%%%%%%%%%%%%%%%%%%%%%%%%%%%%%%%%%%%%%%%%%%%%%%%%%%%%%%%%%%%%%%%%%%%%%%%%%%%%%%%%%%%%%%%%%%

\vspace{0.1cm} 
The group $ \mathrm{SU}((n,1),\C)$ acts on $\mathbb{B}^n$ via bi-holomorphic mappings, which are given by fractional linear transformations, which we explain here. Let 
\begin{align*}
\gamma=\left(\begin{array}{cc} A&B\\C &D\end{array}\right)\in \mathrm{SU}((n,1),\C),
\end{align*}
%%%%%%%%%%%%%%%%%%%%%%%%%%%%%%%%%%%%%%%%%%%%%%%%%%%%%%%%%%%%%%%%%%%%%%%%%%%%%%%%%%%%%%%%%%%
where $A\in M_{n\times n}(\C)$, $B\in M_{n\times 1}(\C)$, $C\in M_{1\times n}(\C) $, and $D\in\C$. For $z=(z_1,\ldots,z_n)^{t}\in\Bn$, the action of $\gamma$ on $\mathbb{B}^n$, is given by
\begin{align}\label{gammact}
\gamma z:=\frac{Az+B}{Cz+D}.
\end{align}
%%%%%%%%%%%%%%%%%%%%%%%%%%%%%%%%%%%%%%%%%%%%%%%%%%%%%%%%%%%%%%%%%%%%%%%%%%%%%%%%%%%%%%%%%%%
Here, $z$ is seen as a column vector in $\C^n$. 
%%%%%%%%%%%%%%%%%%%%%%%%%%%%%%%%%%%%%%%%%%%%%%%%%%%%%%%%%%%%%%%%%%%%%%%%%%%%%%%%%%%%%%%%%%%

\vspace{0.1cm} 
Let $\Gamma\subset \mathrm{SU}((n,1),\calo_{K})$ be a torsion-free, finite-index subgroup, with only one cusp at $\infty:=(1,0,\ldots,0)^{t}\in \Bn$, where $K$ is totally imaginary number field of degree $2$. Let $X_{\G}:=\Gamma\backslash\mathbb{B}^n$ denote the quotient space, which is called a Picard variety of dimension $n$. The Picard variety $X_{\G}$ admits the structure of a noncompact complex manifold of dimension $n$, with a cuspidal singularity ay $\infty$. The hyperbolic metric $\hypbn(z)$ is $\mathrm{SU}((n,1),\C)$-invariant, and hence, descends to define a K\"ahler metric on $X_{\G}$. Furthermore, the Picard variety $X_{\G}$ is of finite hyperbolic volume.
%%%%%%%%%%%%%%%%%%%%%%%%%%%%%%%%%%%%%%%%%%%%%%%%%%%%%%%%%%%%%%%%%%%%%%%%%%%%%%%%%%%%%%%%%%%

\vspace{0.1cm}
Locally, we identify the Picard variety $X_{\G}$ with its universal cover $\Bn$. Locally, for $z,w\in X_{\G}$ (identifying $X_{\G}$ with $\Bn$), the hyperbolic distance between the points $z$ and $w$, is given by formula \eqref{coshreln}. Furthermore, let $\infty=(1,0,\ldots,0)^{t}\in\calf_{\G}$ denote a fixed fundamental domain of $X_{\G}$, and let $\overline{\calf}_{\G}:=\calf_{\G}\cup\lbrace \infty\rbrace$.
%%%%%%%%%%%%%%%%%%%%%%%%%%%%%%%%%%%%%%%%%%%%%%%%%%%%%%%%%%%%%%%%%%%%%%%%%%%%%%%%%%%%%%%%%%%

\vspace{0.1cm}
Set $\G_{3}:=\gamma_{13}\G\gamma_{31}$, then the quotient space $X_{\G_3}:=\G_{3}\backslash \E^n$ is isometric to $X_{\G}$, as a complex manifold. Let $\calf_{\G_{3}}$ denote a fixed fundamental domain of $X_{\G_{3}}$, and recall that $\gamma_{31}^{-1}=\gamma_{13}$, and $\calf_{\G_{3}}=\gamma_{13}(\calf_{\G})$. For the rest of the article, we identify $X_{\G_{3}}$ with $X_{\G}$.
%%%%%%%%%%%%%%%%%%%%%%%%%%%%%%%%%%%%%%%%%%%%%%%%%%%%%%%%%%%%%%%%%%%%%%%%%%%%%%%%%%%%%%%%%%%
%%%%%%%%%%%%%%%%%%%%%%%%%%%%%%%%%%%%%%%%%%%%%%%%%%%%%%%%%%%%%%%%%%%%%%%%%%%%%%%%%%%%%%%%%%%

\vspace{0.2cm} 
\subsection{Stabilizers of the cusps and cuspidal neighborhood of $X_{\G}$}\label{subsec-2.4}
%%%%%%%%%%%%%%%%%%%%%%%%%%%%%%%%%%%%%%%%%%%%%%%%%%%%%%%%%%%%%%%%%%%%%%%%%%%%%%%%%%%%%%%%%%%
We now describe $\G_{3,\infty}$, the isotropy subgroup of $\G_3$, which fixes $\infty=(-\infty,0,\ldots,0)^{t}\in \E^n$. We now consider model (3) of the complex hyperbolic space. The boundary $\partial \overline{\E^n}$ of $\overline{\E^n}$ is $\mathbb P (V_0 (H_3))$ with a distinguished point $\infty$. The Bailey--Borel compactification of $X_{\G}$ is obtained by adding orbits of 
rational boundary points, which are known as the cusps. As in the classical case, there are isometries acting transitively on the cusps. Hence, we only consider the cusp $\infty$. 
%%%%%%%%%%%%%%%%%%%%%%%%%%%%%%%%%%%%%%%%%%%%%%%%%%%%%%%%%%%%%%%%%%%%%%%%%%%%%%%%%%%%%%%%%%%

\vspace{0.1cm} 
Let $\mathrm{U}((n,1),\C)_\infty$ denote the stabilizer of the cusp $\infty$ in $\mathrm{U}((n,1),\C)$. Then any $P \in \mathrm{U}\big((n,1),\C\big)_\infty$ can be decomposed as a product of Heisenberg translation, dilation and rotation. Then 
\begin{align*}
P = \begin{pmatrix} r & -\tau^* U & \frac{-|\tau|^2 + it}{2r} \\ 0 & U & \tau/r \\ 0 & 0 & 1/r  \end{pmatrix},
\end{align*}
%%%%%%%%%%%%%%%%%%%%%%%%%%%%%%%%%%%%%%%%%%%%%%%%%%%%%%%%%%%%%%%%%%%%%%%%%%%%%%%%%%%%%%%%%%%
where $r \in \R^+$, $U \in U(n-1)$, $\tau \in \C^{n-1}$ and $t \in \R$. In particular, for a discrete subgroup $\G_{3} \subset \mathrm{U}\big((n,1),\C\big)$ the group
\begin{align*}
\G_{3,\infty}: = \mathrm{U}((n,1),\C)_\infty \cap \G_{3}
\end{align*}
%%%%%%%%%%%%%%%%%%%%%%%%%%%%%%%%%%%%%%%%%%%%%%%%%%%%%%%%%%%%%%%%%%%%%%%%%%%%%%%%%%%%%%%%%%%
is the stabilizer of the cusp $\infty$ in $\G_{3}$.
%%%%%%%%%%%%%%%%%%%%%%%%%%%%%%%%%%%%%%%%%%%%%%%%%%%%%%%%%%%%%%%%%%%%%%%%%%%%%%%%%%%%%%%%%%%

\vspace{0.1cm} 
Hence, we conclude that $\G_{3,\infty}$ 
\begin{align}\label{gamma3infty}
\G_{3,\infty} \subset \left\lbrace \begin{pmatrix} 1 & -\tau^* U & \frac{-|\tau|^2 + it}{2} \\ 0 & U & \tau \\ 0 & 0 & 1 \end{pmatrix}\Bigg|\,  
U \in \mathrm{SU} (n-1,\calo_K),\,\,\tau \in \calo_K^{n-1} ,\,\,t \in \R  \right\rbrace
\end{align} 
%%%%%%%%%%%%%%%%%%%%%%%%%%%%%%%%%%%%%%%%%%%%%%%%%%%%%%%%%%%%%%%%%%%%%%%%%%%%%%%%%%%%%%%%%%%
is a finite index subgroup.
%%%%%%%%%%%%%%%%%%%%%%%%%%%%%%%%%%%%%%%%%%%%%%%%%%%%%%%%%%%%%%%%%%%%%%%%%%%%%%%%%%%%%%%%%%%

\vspace{0.1cm}
\begin{rem}\label{rem-2.1}
Let us make the group $\G_{3,\infty}$ precise when $\G_3 = \mathrm{SU} ((n,1), \calo_K)$. We have $\tau \in \calo_K^{n-1}$ and $U \in \mathrm{SU} (n-1,\calo_K)$ as before. We only need to determine a more precise condition on $t$. Since $\tau \in \calo_K^{n-1}$, we have $|\tau|^2 \in \calo_K \cap \R_+ = \Z_+$.  Therefore, the condition $\frac{-|\tau|^2 + it}{2} \in \calo_K$ is equivalent to $it \in \calo_K$ and $|\tau|^4 + t^2 \in 4\Z$. 
%%%%%%%%%%%%%%%%%%%%%%%%%%%%%%%%%%%%%%%%%%%%%%%%%%%%%%%%%%%%%%%%%%%%%%%%%%%%%%%%%%%%%%%%%%%

Now, if we take $K=\Q[\sqrt{-d}]$, then we have $\calo_K = \Z [\sqrt{-d}]$ if $d \equiv 1 \textrm{ or } 2 \pmod 4$ and $\calo_K = \Z [\frac{1+\sqrt{-d}}{2}]$ if $d \equiv 3 \pmod 4$. In either case, the condition $it \in \calo_K$ is equivalent to $t = b \sqrt{d}$ for some $b \in \Z$.  When $d \equiv 3 \pmod 4$, the second condition implies that the set of all possible $\tau$ is full lattice $\calo_K^{n-1}$ and all possible $t$ is $\Z \sqrt{d}$. When $d \equiv 1 \textrm{ or } 2 \pmod 4$, the second condition implies that the set of all possible $\tau$ is  the sublattice of index $2$ in $\calo_K^{n-1}$ given by $\{\tau \in \calo_K^{n-1} \mid |\tau|^2 \in 2 \Z \}$ and all possible $t$ is $2\Z \sqrt{d}$.
\end{rem}
%%%%%%%%%%%%%%%%%%%%%%%%%%%%%%%%%%%%%%%%%%%%%%%%%%%%%%%%%%%%%%%%%%%%%%%%%%%%%%%%%%%%%%%%%%%

\vspace{0.1cm} 
The point $\infty=(-\infty,0,\ldots,0)^{t}$ corresponds to a cuspidal singularity for the complex manifold $X_{\G_{3}}$, which is identified with its fundamental domain  $\calf_{\G_{3}}$, a  neighbordhood around $\infty$ is described by the following equation 
\begin{align}
&U_{\epsilon}(\infty):=\big\lbrace z:=(z_1,\ldots,z_n)^{t}\in \calf_{\G_{3}}\big|\,2\mathrm{Re}(z_1)\leq-1\slash \epsilon,\,\,z_2=\cdots=z_n=0 \big\rbrace;\notag\\
&\overline{U}_{\epsilon}(\infty):=U_{\epsilon}(\infty)\cup\lbrace\infty\rbrace.\label{cusp-nbd}
\end{align}
%%%%%%%%%%%%%%%%%%%%%%%%%%%%%%%%%%%%%%%%%%%%%%%%%%%%%%%%%%%%%%%%%%%%%%%%%%%%%%%%%%%%%%%%%%%
%%%%%%%%%%%%%%%%%%%%%%%%%%%%%%%%%%%%%%%%%%%%%%%%%%%%%%%%%%%%%%%%%%%%%%%%%%%%%%%%%%%%%%%%%%%

\vspace{0.2cm}
\subsection{Determinant bundle of the Cotangent bundle}
%%%%%%%%%%%%%%%%%%%%%%%%%%%%%%%%%%%%%%%%%%%%%%%%%%%%%%%%%%%%%%%%%%%%%%%%%%%%%%%%%%%%%%%%%%%
Let $\overline{X}_{\G}:=X_{\G}\cup\lbrace \infty\rbrace$ denote the one-point compactification of $X_{\G}$. Let $\Lambda:=\Omega^n_{\overline{X}_{\G}}$ denote the line bundle of holomorphic $(n,0)$-forms on $\overline{X}_{\G}$. Furthermore, let $\lambda^{k}:=\Lambda^{\otimes k}\otimes O_{\overline{X}_{\G}}((k-1)\infty)$, denote the space of meromorphic differential forms, with a pole of order $k-1$ at the cusp $\infty$, and which remain holomorphic on $\overline{X}_{\G}\backslash\lbrace \infty\rbrace$. 
%%%%%%%%%%%%%%%%%%%%%%%%%%%%%%%%%%%%%%%%%%%%%%%%%%%%%%%%%%%%%%%%%%%%%%%%%%%%%%%%%%%%%%%%%%%

For any $k\geq 1$, let $H^{0}(\overline{X}_{\G},\lambda^{k})$ denote the complex vector space of global holomorphic sections of the line bundle $\lambda^{k}$.  The hyperbolic metric induces a point-wise metric on $H^{0}(\overline{X}_{\G},\lambda^{k})$, which we denote by 
$| \cdot|_{\mathrm{hyp}}$, and an $L^2$-metric, which for any 
\begin{align*}
\omega_1:=f_1(z)\big(dz_1\wedge\cdots\wedge dz_n\big)^{\otimes k},\,\,\,\omega_2:=f_2(z)\big(dz_1\wedge\cdots\wedge dz_n\big)^{\otimes k}\in H^{0}(X_{\G},\lambda^{ k}), 
\end{align*}
%%%%%%%%%%%%%%%%%%%%%%%%%%%%%%%%%%%%%%%%%%%%%%%%%%%%%%%%%%%%%%%%%%%%%%%%%%%%%%%%%%%%%%%%%%%
is given the following equation
\begin{align*}
\langle \omega_1,\omega_2\rangle_{\mathrm{hyp}}:=\int_{\calf_{\G}}(1-|z|^2)^{k(n+1)}f_1(z)\overline{f_2(z)}\hypbnvol(z).
\end{align*} 
%%%%%%%%%%%%%%%%%%%%%%%%%%%%%%%%%%%%%%%%%%%%%%%%%%%%%%%%%%%%%%%%%%%%%%%%%%%%%%%%%%%%%%%%%%%

Let $d_{k}$ denote the dimension of $H^{0}(\overline{X}_{\G},\lambda^{k})$, with respect to the $L^2$ inner-product $\langle \cdot,\cdot\rangle_{\mathrm{hyp}}$, and let $\lbrace \omega_1,\ldots\omega_{d_{k}}\rbrace$ denote an orthonormal basis of $H^{0}(\overline{X}_{\G},\lambda^{k})$. From the Riemann-Roch theorem, it is clear that $d_{k}$ satisfies the following estimate
\begin{align}\label{dk}
d_{k}=O_{X_{\G}}(k^n),
\end{align}
%%%%%%%%%%%%%%%%%%%%%%%%%%%%%%%%%%%%%%%%%%%%%%%%%%%%%%%%%%%%%%%%%%%%%%%%%%%%%%%%%%%%%%%%%%%
where the implied constant depends only on $X_{\G}$.
%%%%%%%%%%%%%%%%%%%%%%%%%%%%%%%%%%%%%%%%%%%%%%%%%%%%%%%%%%%%%%%%%%%%%%%%%%%%%%%%%%%%%%%%%%%

\vspace{0.1cm}
The Bergman kernel associated to the complex vector space $H^{0}(\overline{X}_{\G},\lambda^{k})$ is given by the following formula
\begin{align*}
\bkx(z,w):=\sum_{j=1}^{d_{k}}\omega_{j}(z)\overline{\omega_{j}(w)}.
\end{align*}
%%%%%%%%%%%%%%%%%%%%%%%%%%%%%%%%%%%%%%%%%%%%%%%%%%%%%%%%%%%%%%%%%%%%%%%%%%%%%%%%%%%%%%%%%%%
From Riesz representation, it follows that the definition of the Bergman kernel $\bkx$, is independent of the choice of orthonormal basis of  $H^{0}(\overline{X}_{\G},\lambda^{ k})$. 
%%%%%%%%%%%%%%%%%%%%%%%%%%%%%%%%%%%%%%%%%%%%%%%%%%%%%%%%%%%%%%%%%%%%%%%%%%%%%%%%%%%%%%%%%%%

\vspace{0.1cm}
The Bergman kernel $\bkx$ is a holomorphic $k$-differential form in the $z$-variable on $X_{\G}$, and an antiholomorphic $k$-differential form in the $w$-variable on $X_{\G}$. Furthermore, the hyperbolic metric $| \cdot|_{\mathrm{hyp}}$ of the Bergman kernel, is given by the following formula
\begin{align*}
\big|\bkx(z,w) \big|_{\mathrm{hyp}}:=(1-|z|^{2})^{k(n+1)\slash 2}(1-|w|^{2})^{k(n+1)\slash 2}\big|\bkx(z,w) \big|.
\end{align*} 
%%%%%%%%%%%%%%%%%%%%%%%%%%%%%%%%%%%%%%%%%%%%%%%%%%%%%%%%%%%%%%%%%%%%%%%%%%%%%%%%%%%%%%%%%%%

\vspace{0.1cm}
From the isometry $X_{\G}\simeq X_{\G_{3}}$, and identifying $X_{\G}$ with $X_{\G_3}$, for $z,w \in X_{\G_3}$, it follows that 
 \begin{align*}
\big|\bkx(z,w) \big|_{\mathrm{hyp}}:=\big(-\mathrm{Re}(z_1)-\sum_{j=2}^{n}|z_{j}|^{2}\big)^{k(n+1)\slash 2}\big(-\mathrm{Re}(w_1)-\sum_{j=2}^{n}|w_{j}|^{2}\big)^{k(n+1)\slash 2}\big|\bkx(z,w) \big|.
\end{align*} 
%%%%%%%%%%%%%%%%%%%%%%%%%%%%%%%%%%%%%%%%%%%%%%%%%%%%%%%%%%%%%%%%%%%%%%%%%%%%%%%%%%%%%%%%%%%
%%%%%%%%%%%%%%%%%%%%%%%%%%%%%%%%%%%%%%%%%%%%%%%%%%%%%%%%%%%%%%%%%%%%%%%%%%%%%%%%%%%%%%%%%%%

\vspace{0.2cm}
\subsection{Cusp forms and Bergman kernel}
%%%%%%%%%%%%%%%%%%%%%%%%%%%%%%%%%%%%%%%%%%%%%%%%%%%%%%%%%%%%%%%%%%%%%%%%%%%%%%%%%%%%%%%%%%%
For any $n\geq 2$, a holomorphic function $f:\Bn\longrightarrow \C$ is called a cusp form of weight $k\geq 1$, with respect to $\G$, if $f$ satisfies the following transformation property
\begin{align*}
&f(\gamma z)=(CZ+D)^{k}f(z), \,\,\mathrm{where}\,\,\gamma=\left(\begin{smallmatrix}A &B\\C&D  \end{smallmatrix}\right)\in\G,
\end{align*}
%%%%%%%%%%%%%%%%%%%%%%%%%%%%%%%%%%%%%%%%%%%%%%%%%%%%%%%%%%%%%%%%%%%%%%%%%%%%%%%%%%%%%%%%%%%
where $A\in M_{n\times n}(\C)$, $B\in M_{n\times 1}(\C)$, $C\in M_{1\times n}(\C)$, $D\in \C$,  and $f(\infty)=0$.
%%%%%%%%%%%%%%%%%%%%%%%%%%%%%%%%%%%%%%%%%%%%%%%%%%%%%%%%%%%%%%%%%%%%%%%%%%%%%%%%%%%%%%%%%%%

\vspace{0.1cm}
Any $f\in\sk$, from the Fourier-Jacobi expansion of Picard cusp forms stated in section 4 of \cite{geer}, after adjusting the normalizations, we have
\begin{align}\label{fourier}
f(\gamma_{31}z)=\sum_{n=1}^{\infty}f_{n}(z_2,\ldots,z_n)e^{2\pi c_{\G} nz_{1}},
\end{align}
%%%%%%%%%%%%%%%%%%%%%%%%%%%%%%%%%%%%%%%%%%%%%%%%%%%%%%%%%%%%%%%%%%%%%%%%%%%%%%%%%%%%%%%%%%%
where $\gamma_{31}$ is as defined in \eqref{gamma31}, and $z:=(z_{1},\ldots,z_{n})\in\E^n$, and the functions $\lbrace f_n\rbrace_{n\geq 1}$ are holomorphic functions, and $c_{\G}$ is a constant which depends only on $\G$. Without loss of generality, for the brevity of ensuing computations, we choose $c_{\G} = 1$.
%%%%%%%%%%%%%%%%%%%%%%%%%%%%%%%%%%%%%%%%%%%%%%%%%%%%%%%%%%%%%%%%%%%%%%%%%%%%%%%%%%%%%%%%%%%

\vspace{0.1cm}
For any $k\geq 1$, let $\sk$ denote the complex vector space of weight-$k$ cusp forms with respect to $\G$. For any $f,g\in\sk$ the Petersson inner-product is given by the 
following formula
\begin{align}\label{pet-ip}
\big\langle f,g\big\rangle_{\mathrm{pet}}:=\int_{\calf_{\G}}(1-|z|^2)^{k}f(z)\overline{g(z)}\hypbnvol(z),
\end{align}
%%%%%%%%%%%%%%%%%%%%%%%%%%%%%%%%%%%%%%%%%%%%%%%%%%%%%%%%%%%%%%%%%%%%%%%%%%%%%%%%%%%%%%%%%%%
where $\calf_{\G}$ is a fixed fundamental domain of $X_{\G}$. 
%%%%%%%%%%%%%%%%%%%%%%%%%%%%%%%%%%%%%%%%%%%%%%%%%%%%%%%%%%%%%%%%%%%%%%%%%%%%%%%%%%%%%%%%%%%

\vspace{0.1cm}
For any $f\in\sk$, the Petersson norm at the point $z\in X_{\G}$ (identifying $X_{\G}$ with its universal cover $\Bn$), is given by the following formula
\begin{align*}
\big| f(z)\big|^{2}_{\mathrm{pet}}=\big|\langle \tilde{z},\tilde{z}\rangle_{H}\big|^{k}|f(z)|^{2}=(1-|z|^{2})^{k}|f(z)|^{2}.
\end{align*}
%%%%%%%%%%%%%%%%%%%%%%%%%%%%%%%%%%%%%%%%%%%%%%%%%%%%%%%%%%%%%%%%%%%%%%%%%%%%%%%%%%%%%%%%%%%

\vspace{0.1cm}
For any $k\geq 1$, we have the following isometry 
\begin{align}\label{iso}
H^{0}\big(\overline{X}_{\G},\lambda^{k}\big)=\mathcal{S}_{k(n+1)}(\G).
\end{align}
%%%%%%%%%%%%%%%%%%%%%%%%%%%%%%%%%%%%%%%%%%%%%%%%%%%%%%%%%%%%%%%%%%%%%%%%%%%%%%%%%%%%%%%%%%%

\vspace{0.1cm}
From the above isometry, it follows that the dimension of $\skn$ as a complex vector space is equal to $d_{k}$, where $d_k$ is the dimension of $H^{0}(\overline{X}_{\G},\lambda^{k})$, which satisfies estimate \eqref{dk}.
%%%%%%%%%%%%%%%%%%%%%%%%%%%%%%%%%%%%%%%%%%%%%%%%%%%%%%%%%%%%%%%%%%%%%%%%%%%%%%%%%%%%%%%%%%%

\vspace{0.1cm}
We now describe the Bergman kernel $\bkn$ associated to the complex vector space $\skn$. Let $\lbrace f_{1},\ldots,f_{d_{k}}\rbrace$ denote an orthonormal basis of $\sk$, with respect to the Petersson inner-product, which is as defined in \eqref{pet-ip}. Let $\big\lbrace f_1,\ldots, f_{d_{k}}\big\rbrace$ denote an orthonormal basis of $\skn$. Then, the Bergman kernel associated to $\skn$, is given by the following formula
\begin{align}\label{bkn}
\bkn(z,w):=\sum_{j=1}^{d_{k}}f_{j}(z)\overline{f_{j}(w)}.
\end{align}
%%%%%%%%%%%%%%%%%%%%%%%%%%%%%%%%%%%%%%%%%%%%%%%%%%%%%%%%%%%%%%%%%%%%%%%%%%%%%%%%%%%%%%%%%%%
From Riesz representation, it follows that the definition of the Bergman kernel $\bkn$, is independent of the choice of orthonormal basis of  $\skn$. 
%%%%%%%%%%%%%%%%%%%%%%%%%%%%%%%%%%%%%%%%%%%%%%%%%%%%%%%%%%%%%%%%%%%%%%%%%%%%%%%%%%%%%%%%%%%

\vspace{0.1cm}
The Bergman kernel $\bkn$ is a holomorphic cusp form of weight-$k$ in the $z$-variable, and an antiholomorphic cusp form of weight-$k$ in the $w$-variable. Furthermore, the Petersson metric of the Bergman kernel, is given by the following formula
\begin{align}\label{bkhyp}
\big|\bkn(z,w) \big|_{\mathrm{hyp}}:=\big(1-|z|^{2}\big)^{k(n+1)\slash 2}\big(1-|w|^{2}\big)^{k(n+1)\slash 2}\big|\bkn(z,w) \big|.
\end{align} 
%%%%%%%%%%%%%%%%%%%%%%%%%%%%%%%%%%%%%%%%%%%%%%%%%%%%%%%%%%%%%%%%%%%%%%%%%%%%%%%%%%%%%%%%%%%

\vspace{0.1cm}
From isometry \eqref{iso}, for any $k\geq (n+1)$  it follows that 
\begin{align}\label{bkiso}
\big| \bkx(z,w)\big|_{\mathrm{hyp}}=\big| \bkn(z,w)\big|_{\mathrm{pet}}.
\end{align}
%%%%%%%%%%%%%%%%%%%%%%%%%%%%%%%%%%%%%%%%%%%%%%%%%%%%%%%%%%%%%%%%%%%%%%%%%%%%%%%%%%%%%%%%%%%

\vspace{0.1cm}
Moreover, for any $z,w\in \Bn$ the Bergman kernel $\bkn$ admits the following infinite series expansion
\begin{align}\label{bkseries}
\bkn(z,w):=\sum_{\gamma\in \G}\frac{\ckn}{\big(\langle \tilde{z},\tilde{ \gamma w}\rangle_{H}\big(\overline{Cw+D}\big)\big)^{k(n+1)}}, \,\,\,\,\mathrm{where}\,\,\gamma=\left(\begin{smallmatrix}A &B\\C&D  \end{smallmatrix}\right)\in\G,
\end{align}
%%%%%%%%%%%%%%%%%%%%%%%%%%%%%%%%%%%%%%%%%%%%%%%%%%%%%%%%%%%%%%%%%%%%%%%%%%%%%%%%%%%%%%%%%%%
where $A\in M_{n\times n}(\C)$, $B\in M_{n\times 1}(\C)$, $C\in M_{1\times n}(\C)$, $D\in \C$, and $\tilde{z}, \tilde{\gamma w}$ are lifts of $z$ and $\gamma w$ to $\cala^{n+1}$, and the inner-product $\langle \cdot,\cdot\rangle_{H}$ is as defined in equation \eqref{coshreln}. The constant $\ckn$ depends only on $X_{\G}$ and $k$. In Proposition \ref{prop2}, we will show that the constant $\ckn$ satisfies the following estimate
\begin{align}\label{ck}
\ckn=O_{X_{\G}}(k^{n}),
\end{align}
%%%%%%%%%%%%%%%%%%%%%%%%%%%%%%%%%%%%%%%%%%%%%%%%%%%%%%%%%%%%%%%%%%%%%%%%%%%%%%%%%%%%%%%%%%%
where the implied constant depends only on $X_{\G}$.
%%%%%%%%%%%%%%%%%%%%%%%%%%%%%%%%%%%%%%%%%%%%%%%%%%%%%%%%%%%%%%%%%%%%%%%%%%%%%%%%%%%%%%%%%%%

\vspace{0.1cm}
Furthermore, combining equations \eqref{bkseries} and 
\eqref{bkhyp}, we arrive at the following estimate
\begin{align*}
\big|\bkn(z,w)\big|_{\mathrm{pet}}\leq \sum_{\gamma\in \G}\frac{\ckn(1-|z|^{2})^{k(n+1)\slash 2}(1-|w|^{2})^{k(n+1)\slash 2}}{\big|\langle \tilde{z},\tilde{ \gamma w}\rangle_{H}\big|^{k(n+1)}\cdot\big|Cw+D\big|^{k(n+1)}}.
\end{align*}
%%%%%%%%%%%%%%%%%%%%%%%%%%%%%%%%%%%%%%%%%%%%%%%%%%%%%%%%%%%%%%%%%%%%%%%%%%%%%%%%%%%%%%%%%%%
Combining the above inequality with inequality (24) in \cite{anil7}, we infer that 
\begin{align}\label{bkineq}
\big|\bkn(z,w)\big|_{\mathrm{pet}}\leq  \sum_{\gamma\in \G}\frac{\ckn}{\cosh^{k(n+1)}(\dhyp(z,\gamma w)\slash 2)}.
\end{align} 
%%%%%%%%%%%%%%%%%%%%%%%%%%%%%%%%%%%%%%%%%%%%%%%%%%%%%%%%%%%%%%%%%%%%%%%%%%%%%%%%%%%%%%%%%%%
%%%%%%%%%%%%%%%%%%%%%%%%%%%%%%%%%%%%%%%%%%%%%%%%%%%%%%%%%%%%%%%%%%%%%%%%%%%%%%%%%%%%%%%%%%%

\vspace{0.2cm}
\subsection{Bergman metric}
%%%%%%%%%%%%%%%%%%%%%%%%%%%%%%%%%%%%%%%%%%%%%%%%%%%%%%%%%%%%%%%%%%%%%%%%%%%%%%%%%%%%%%%%%%%
For any $k\geq 1$, and $z\in X_{\G}$, the Bergman metric associated to the line bundle $\lambda^k$, is given by the following formula
\begin{align}\label{bergman-metric-defn}
\mu_{\mathrm{Ber},k}(z):=-\frac{i}{2\pi}\partial_{z}\partial_{\overline{z}}\big|\bkx(z,z)\big|_{\mathrm{pet}},
\end{align}
%%%%%%%%%%%%%%%%%%%%%%%%%%%%%%%%%%%%%%%%%%%%%%%%%%%%%%%%%%%%%%%%%%%%%%%%%%%%%%%%%%%%%%%%%%%
and the associated volume form is given by the following formula
\begin{align*}
\mu_{\mathrm{Ber},k}^{\mathrm{vol}}(z):=\frac{1}{n!}\bigwedge_{j=1}^{n}\mu_{\mathrm{Ber},k}(z).
\end{align*}
%%%%%%%%%%%%%%%%%%%%%%%%%%%%%%%%%%%%%%%%%%%%%%%%%%%%%%%%%%%%%%%%%%%%%%%%%%%%%%%%%%%%%%%%%%%

Furthermore, from equation (28) from \cite{anil7}, we arrive at the following equation
\begin{align}\label{bergman-metric-1}
\mu_{\mathrm{Ber},k}(z)=k(n+1)\hyp(z)-\frac{i}{2\pi}\partial_{z}\partial_{\overline{z}}\big|\bkn(z,z) \big|.
\end{align}
%%%%%%%%%%%%%%%%%%%%%%%%%%%%%%%%%%%%%%%%%%%%%%%%%%%%%%%%%%%%%%%%%%%%%%%%%%%%%%%%%%%%%%%%%%%
%%%%%%%%%%%%%%%%%%%%%%%%%%%%%%%%%%%%%%%%%%%%%%%%%%%%%%%%%%%%%%%%%%%%%%%%%%%%%%%%%%%%%%%%%%%

\vspace{0.2cm}
\subsection{Counting functions of $\G$ and auxiliary estimates from complex geometry}
%%%%%%%%%%%%%%%%%%%%%%%%%%%%%%%%%%%%%%%%%%%%%%%%%%%%%%%%%%%%%%%%%%%%%%%%%%%%%%%%%%%%%%%%%%%
We will now describe a counting function, which we will be used to estimate elements of $\G$ in section \ref{sec-3}. Recall that $(1,0,\ldots,0)^{t}\calf_{\G}$  denotes a fixed fundamental 
domain of $X_{\G}$. 
%%%%%%%%%%%%%%%%%%%%%%%%%%%%%%%%%%%%%%%%%%%%%%%%%%%%%%%%%%%%%%%%%%%%%%%%%%%%%%%%%%%%%%%%%%%

\vspace{0.1cm}
For any $z,w\in \calf_{\G}$, set
\begin{align*}
&\cals_{\G}(z,w,;\delta):=\big\lbrace \gamma\in\G\backslash \G_{\infty}\big|\,\dhyp(z,\gamma w\big)\leq \delta\big\rbrace;\\[0.12cm]
&\caln_{\G}(z,w;\delta):=\mathrm{card}\,\cals_{\G}(z,w,;\delta).
\end{align*}
%%%%%%%%%%%%%%%%%%%%%%%%%%%%%%%%%%%%%%%%%%%%%%%%%%%%%%%%%%%%%%%%%%%%%%%%%%%%%%%%%%%%%%%%%%%

\vspace{0.1cm}
Furthermore, the injectivity radius of $X_{\G}$ is given by the following formula
\begin{align}\label{irad}
&\rx:=\inf\big\lbrace\dhyp(z,\gamma z)\big|\,z\in \calf_{\G},\,\,\gamma\in\G\backslash \G_{\infty} \big\rbrace.
\end{align}
%%%%%%%%%%%%%%%%%%%%%%%%%%%%%%%%%%%%%%%%%%%%%%%%%%%%%%%%%%%%%%%%%%%%%%%%%%%%%%%%%%%%%%%%%%%

\vspace{0.1cm}
We denote the volume of a geodesic ball of radius $r>0$, around a given $z\in \Bn$, by $\mathrm{vol}(B(z,r))$, and it is given by the following formula
\begin{align*}
\mathrm{vol}(B(z,r))=\frac{(4\pi)^{n}}{n!}\sinh^{2n}(r\slash 2).
\end{align*}  
%%%%%%%%%%%%%%%%%%%%%%%%%%%%%%%%%%%%%%%%%%%%%%%%%%%%%%%%%%%%%%%%%%%%%%%%%%%%%%%%%%%%%%%%%%%

\vspace{0.1cm}
We now state estimates of $\caln_{\G}$ from Lemma 4.1 in \cite{anil4}, after replacing some of the constants appearing in the estimates by their optimal versions. For any 
$z,w\in \calf_{\G}$, and $\delta>0$, we have the following estimates
\begin{align}\label{cfn1}
\caln_{\G}(z,w;\delta)\leq \frac{\sinh^{2n}((2\delta+\rx)\slash 4)}{\sinh^{2n}(\rx\slash 4)}.
\end{align}
%%%%%%%%%%%%%%%%%%%%%%%%%%%%%%%%%%%%%%%%%%%%%%%%%%%%%%%%%%%%%%%%%%%%%%%%%%%%%%%%%%%%%%%%%%%

Furthermore, let $f$ be any positive, smooth, monotonically decreasing function defined on $\R_{\geq 0}$.  Then, for any $\delta>\rx\slash 2$, assuming that all the integral exists, 
we have the following estimate
\begin{align}\label{cfn2}
\int_{0}^{\infty}f(\rho)d\caln_{\G}(,w;\rho)\leq \int_{0}^{\delta}f(\rho)d\caln_{\G}(z,w;\rho)+f(\delta) \frac{\sinh^{2n}((2\delta+\rx)\slash 4)}{\sinh^{2n}(\rx\slash 4)}+\notag\\[0.12cm]
\frac{2n}{\sinh^{2n}(\rx\slash 4)}\int_{\delta}^{\infty}f(\rho)\sinh^{2n-1}((2\rho+\rx\slash 4)\slash 4)\cosh((2\rho+\rx)\slash4).
\end{align}
%%%%%%%%%%%%%%%%%%%%%%%%%%%%%%%%%%%%%%%%%%%%%%%%%%%%%%%%%%%%%%%%%%%%%%%%%%%%%%%%%%%%%%%%%%%

\vspace{0.1cm}
From complex geometry, we now state an estimate of the Bergman kernel, associated to a holomorphic line bundle, defined over a complex manifold. Let $X$ be a nocompact 
complex manifold, and let $\lambda$ and $\mathcal{W}$ be a holomorphic line bundle and a holomorphi vector bundle of rank $r$, respectively, defined over $X$. For any $k\geq 1$, let $\mathcal{B}_{X}^{\lambda^{k}, \mathcal{W}}$ denote the Bergman kernel associated to $H^{0}(X,\lambda^{\otimes k}\otimes \mathcal{W})$, the complex vector space of holomorphic global sections of the vector bundle $\lambda^{\otimes k}\otimes \mathcal{W}$. Furthermore, let $|\cdot|_{\lambda^{\otimes k}\otimes \mathcal{W}}$ denote the point-wise metric on $H^{0}(X,\lambda^{\otimes k}\otimes \mathcal{W})$.
%%%%%%%%%%%%%%%%%%%%%%%%%%%%%%%%%%%%%%%%%%%%%%%%%%%%%%%%%%%%%%%%%%%%%%%%%%%%%%%%%%%%%%%%%%%

\vspace{0.1cm}
Then, in \cite{berman}, Berman proved the following estimate 
\begin{align}\label{berman1}
\limsup_{k\rightarrow \infty}\frac{1}{k^{n}}\big|\mathcal{B}_{X}^{\lambda^{k},\mathcal{W}}(z,z)\big|_{\lambda^{\otimes k}\otimes\mathcal{W}}=O_{X}(1),
\end{align}  
%%%%%%%%%%%%%%%%%%%%%%%%%%%%%%%%%%%%%%%%%%%%%%%%%%%%%%%%%%%%%%%%%%%%%%%%%%%%%%%%%%%%%%%%%%%
where the implied constant depends only on $X$. From estimate \eqref{berman1}, for any compact subset $A\subset X$, we have the following estimate
\begin{align}\label{berman2}
\big|\mathcal{B}_{X}^{\lambda^{k},\mathcal{W}}(z,z)\big|_{\lambda^{\otimes k}\otimes\mathcal{W}}=O_{A}\big(r\cdot k^{n}\big),
\end{align}
%%%%%%%%%%%%%%%%%%%%%%%%%%%%%%%%%%%%%%%%%%%%%%%%%%%%%%%%%%%%%%%%%%%%%%%%%%%%%%%%%%%%%%%%%%%
where the implied constant depends only on $A$.
%%%%%%%%%%%%%%%%%%%%%%%%%%%%%%%%%%%%%%%%%%%%%%%%%%%%%%%%%%%%%%%%%%%%%%%%%%%%%%%%%%%%%%%%%%%
%%%%%%%%%%%%%%%%%%%%%%%%%%%%%%%%%%%%%%%%%%%%%%%%%%%%%%%%%%%%%%%%%%%%%%%%%%%%%%%%%%%%%%%%%%%

\vspace{0.2cm}
\section{Estimates of the Bergman kernel}\label{sec-3}
%%%%%%%%%%%%%%%%%%%%%%%%%%%%%%%%%%%%%%%%%%%%%%%%%%%%%%%%%%%%%%%%%%%%%%%%%%%%%%%%%%%%%%%%%%%
In this section, we prove Theorems \ref{mainthm1} and \ref{mainthm2}.
%%%%%%%%%%%%%%%%%%%%%%%%%%%%%%%%%%%%%%%%%%%%%%%%%%%%%%%%%%%%%%%%%%%%%%%%%%%%%%%%%%%%%%%%%%%

\vspace{0.2cm}
\subsection{Estimates of the Bergman kernel}\label{subsec-3.1}
%%%%%%%%%%%%%%%%%%%%%%%%%%%%%%%%%%%%%%%%%%%%%%%%%%%%%%%%%%%%%%%%%%%%%%%%%%%%%%%%%%%%%%%%%%%
In this section, we prove Theorem \ref{mainthm1}.
%%%%%%%%%%%%%%%%%%%%%%%%%%%%%%%%%%%%%%%%%%%%%%%%%%%%%%%%%%%%%%%%%%%%%%%%%%%%%%%%%%%%%%%%%%%

\vspace{0.1cm}
In the following proposition, we prove estimate \eqref{ck}. 
%%%%%%%%%%%%%%%%%%%%%%%%%%%%%%%%%%%%%%%%%%%%%%%%%%%%%%%%%%%%%%%%%%%%%%%%%%%%%%%%%%%%%%%%%%%

\begin{prop}\label{prop1}
With hypothesis as above, let $\cala\subset \calf_{\G}$ be a compact subset. Then, for any $k\geq 2$, and $z\in \cala$, we have the following estimate
\begin{align}\label{prop1:eqn}
\sup_{z\in\cala}\big| \bkx(z,z)\big|_{\mathrm{hyp}}=O_{\cala}\big(\ckn\big),
\end{align}
%%%%%%%%%%%%%%%%%%%%%%%%%%%%%%%%%%%%%%%%%%%%%%%%%%%%%%%%%%%%%%%%%%%%%%%%%%%%%%%%%%%%%%%%%%%
where the implied constant depends only on the compact subset $\cala\subset \calf_{\G}$.
%%%%%%%%%%%%%%%%%%%%%%%%%%%%%%%%%%%%%%%%%%%%%%%%%%%%%%%%%%%%%%%%%%%%%%%%%%%%%%%%%%%%%%%%%%%
\begin{proof}
Let $\cala\subset \calf_{\G}$ be a fixed compact subset. For any $k\geq 1$, and $z\in\cala\subset \calf_{\G}$, combining equation \eqref{bkiso} with estimate \eqref{bkineq}, we have the following estimate
\begin{align}\label{prop1-eqn1}
\big| \bkx(z,z)\big|_{\mathrm{hyp}}\leq \sum_{\gamma\in\G}\frac{\ckn}{\cosh^{k(n+1)}(\dhyp(z,\gamma z)\slash 2)}.
\end{align}  
%%%%%%%%%%%%%%%%%%%%%%%%%%%%%%%%%%%%%%%%%%%%%%%%%%%%%%%%%%%%%%%%%%%%%%%%%%%%%%%%%%%%%%%%%%%

Applying and adatping estimates \eqref{cfn1} and \eqref{cfn2} (refining estimate (50) in \cite{anil4}), and adapting arguments from Theorem 4.2 in \cite{anil4} to the compact subset $\cala\subset \calf_{\G}$, we arrive at the following estimate 
\begin{align}\label{prop1-eqn2}
\sum_{\gamma\in\G}\frac{\ckn}{\cosh^{k(n+1)}(\dhyp(z,\gamma z)\slash 2)}\leq\ckn+\frac{\ckn\sin^{2n}(5r_{\cala}\slash 8)}{\sinh^{2n}(r_{\cala}\slash 4)\cosh^{k(n+1)}(3r_{\cala}\slash 8)}+\notag\\[0.12cm] \frac{\ckn\cosh^{2n}(r_{\cala}\slash 4)}{(k(n+1)-2n-1)\sinh^{2n}(r_{\cala}\slash 4)
\cosh^{k(n+1)-2n-2}(3r_{\cala}/8)},\\[0.12cm]
\mathrm{where}\,\,r_{\cala}:=\big\lbrace \dhyp(z,\gamma z)\big|\,z\in \cala,\,\,\gamma\in\G\backslash\lbrace\mathrm{Id}\rbrace\big\rbrace,\label{ra}
\end{align}
%%%%%%%%%%%%%%%%%%%%%%%%%%%%%%%%%%%%%%%%%%%%%%%%%%%%%%%%%%%%%%%%%%%%%%%%%%%%%%%%%%%%%%%%%%%
and $\mathrm{Id}$ denotes the identity matrix of $\G$. 
%%%%%%%%%%%%%%%%%%%%%%%%%%%%%%%%%%%%%%%%%%%%%%%%%%%%%%%%%%%%%%%%%%%%%%%%%%%%%%%%%%%%%%%%%%%

\vspace{0.1cm}
Combining estimates \eqref{prop1-eqn1} and \eqref{prop1-eqn2}, completes the proof of the proposition.
\end{proof}
\end{prop}
%%%%%%%%%%%%%%%%%%%%%%%%%%%%%%%%%%%%%%%%%%%%%%%%%%%%%%%%%%%%%%%%%%%%%%%%%%%%%%%%%%%%%%%%%%%

\vspace{0.1cm}
In the following proposition, we prove estimate \eqref{ck}, which is already proved in \cite{anil4}. However, for the convenience of the reader, we re-prove estimate \eqref{ck} here.
%%%%%%%%%%%%%%%%%%%%%%%%%%%%%%%%%%%%%%%%%%%%%%%%%%%%%%%%%%%%%%%%%%%%%%%%%%%%%%%%%%%%%%%%%%%

\vspace{0.1cm}
\begin{prop}\label{prop2}
With hypothesis as above, for $k\gg1$, with constant $\ckn$ as defined in equation \eqref{bkseries}, we have the following estimate
\begin{align*}
\ckn=O_{X_{\G}}(k^{n}),
\end{align*} 
%%%%%%%%%%%%%%%%%%%%%%%%%%%%%%%%%%%%%%%%%%%%%%%%%%%%%%%%%%%%%%%%%%%%%%%%%%%%%%%%%%%%%%%%%%%
where the implied constant depends only on $X_{\G}$.
%%%%%%%%%%%%%%%%%%%%%%%%%%%%%%%%%%%%%%%%%%%%%%%%%%%%%%%%%%%%%%%%%%%%%%%%%%%%%%%%%%%%%%%%%%%
\begin{proof}
Let $\cala \subset \calf_{\G}$ be any fixed compact subset. For any $k\gg1$, and $z\in \cala$, combining equation \eqref{bkiso} with estimate \eqref{bkineq}, we have the following estimate
\begin{align}\label{prop2-eqn1}
\big|\bkx(z,z)\big|_{\mathrm{hyp}}\geq \ckn- \sum_{\gamma\in\G\backslash \lbrace \mathrm{Id}\rbrace}\frac{\ckn}{\cosh^{k(n+1)}(\dhyp(z,\gamma z)\slash 2)},
\end{align}
%%%%%%%%%%%%%%%%%%%%%%%%%%%%%%%%%%%%%%%%%%%%%%%%%%%%%%%%%%%%%%%%%%%%%%%%%%%%%%%%%%%%%%%%%%%
where $\mathrm{Id}$ denotes the identity element of $\G$. 
%%%%%%%%%%%%%%%%%%%%%%%%%%%%%%%%%%%%%%%%%%%%%%%%%%%%%%%%%%%%%%%%%%%%%%%%%%%%%%%%%%%%%%%%%%%

\vspace{0.1cm}
From arguments as in Proposition \ref{prop1} (estimate \eqref{prop1-eqn2}), we have the following estimate
\begin{align}\label{prop2-eqn2}
\sum_{\gamma\in\G\backslash \lbrace \mathrm{Id}\rbrace}\frac{\ckn}{\cosh^{k(n+1)}(\dhyp(z,\gamma z)\slash 2)}\leq\frac{\ckn\sin^{2n}(5r_{\cala}\slash 8)}{\sinh^{2n}(r_{\cala}\slash 4)\cosh^{k(n+1)}(3r_{\cala}\slash 8)}+ \notag \\[0.12cm]\frac{\ckn\cosh^{2n}(r_{\cala}\slash 4)}{(k(n+1)-2n-1)\sinh^{2n}(r_{\cala}\slash 4)
\cosh^{k(n+1)-2n-2}(3r_{\cala}/8)},
\end{align}
%%%%%%%%%%%%%%%%%%%%%%%%%%%%%%%%%%%%%%%%%%%%%%%%%%%%%%%%%%%%%%%%%%%%%%%%%%%%%%%%%%%%%%%%%%%
where $r_{\cala}$ is as defined in equation \eqref{ra}.
%%%%%%%%%%%%%%%%%%%%%%%%%%%%%%%%%%%%%%%%%%%%%%%%%%%%%%%%%%%%%%%%%%%%%%%%%%%%%%%%%%%%%%%%%%%

\vspace{0.1cm}
For $k\gg1$, combining estimates \eqref{prop2-eqn1}, \eqref{prop2-eqn2}, \eqref{berman2}, and \eqref{prop1:eqn}, we arrive at the following estimate
\begin{align*}
&\ckn-\frac{\ckn\sin^{2n}(5r_{\cala}\slash 8)}{\sinh^{2n}(r_{\cala}\slash 4)\cosh^{k(n+1)}(3r_{\cala}\slash 8)}-\notag\\&\frac{\ckn\cosh^{2n}(r_{\cala}\slash 4)}{(k(n+1)-2n-1)\sinh^{2n}(r_{\cala}\slash 4)\cosh^{k(n+1)-2n-2}(3r_{\cala}/8)}
\leq\big|\bkx(z,z)\big|_{\mathrm{hyp}}=O_{\cala}\big(k^{n}\big),
\end{align*} 
%%%%%%%%%%%%%%%%%%%%%%%%%%%%%%%%%%%%%%%%%%%%%%%%%%%%%%%%%%%%%%%%%%%%%%%%%%%%%%%%%%%%%%%%%%%
which completes the proof of the proposition.
\end{proof}
\end{prop}
%%%%%%%%%%%%%%%%%%%%%%%%%%%%%%%%%%%%%%%%%%%%%%%%%%%%%%%%%%%%%%%%%%%%%%%%%%%%%%%%%%%%%%%%%%%

\vspace{0.1cm}
In the following lemma, we show that, when restricted to a cuspidal neighborhood, the Petersson norm of the Bergman kernel attains its maximum, on the boundary of the cuspidal neighborhood. 
%%%%%%%%%%%%%%%%%%%%%%%%%%%%%%%%%%%%%%%%%%%%%%%%%%%%%%%%%%%%%%%%%%%%%%%%%%%%%%%%%%%%%%%%%%%

\vspace{0.1cm}
\begin{lem}\label{lem3}
With hypothesis as above, for any $k\geq 2$, we have
\begin{align}\label{lem3:eqn}
\sup_{z\in \overline{U}_{4\pi\slash k(n+1)}(\infty)}\big|\bkx(\gamma_{31}z,\gamma_{31}z) \big|_{\mathrm{hyp}}=\sup_{z\in\partial  U_{4\pi\slash k(n+1)}(\infty)}\big|\bkx(\gamma_{31}z,\gamma_{31}z) \big|_{\mathrm{hyp}},
\end{align}
%%%%%%%%%%%%%%%%%%%%%%%%%%%%%%%%%%%%%%%%%%%%%%%%%%%%%%%%%%%%%%%%%%%%%%%%%%%%%%%%%%%%%%%%%%%
where $\overline{U}_{4\pi\slash k(n+1)}(\infty)$ is as defined in equation \eqref{cusp-nbd}, and $\partial \overline{U}_{4\pi\slash k(n+1)}$ denotes the boundary of $U_{4\pi\slash k(n+1)}(\infty)$.
%%%%%%%%%%%%%%%%%%%%%%%%%%%%%%%%%%%%%%%%%%%%%%%%%%%%%%%%%%%%%%%%%%%%%%%%%%%%%%%%%%%%%%%%%%%
\begin{proof}
From the definition of the Bergman kernel, we have
\begin{align*}
\big|\bkn(\gamma_{31}z,\gamma_{31}z) \big|_{\mathrm{pet}}=\big( -2\mathrm{Re}(z_1)-\sum_{i=2}^{n}|z_{i}|^{2}\big)^{k(n+1)}e^{4\pi \mathrm{Re}(z_1)}\cdot\sum_{j=1}^{d_{k}}\bigg|\frac{f_{j}(z)}{e^{2\pi z_1}}\bigg|^{2}.
\end{align*}
%%%%%%%%%%%%%%%%%%%%%%%%%%%%%%%%%%%%%%%%%%%%%%%%%%%%%%%%%%%%%%%%%%%%%%%%%%%%%%%%%%%%%%%%%%%

From the Fourier expansion of cusp forms, which is as described in equation \eqref{fourier}, we infer that the term
\begin{align*}
\sum_{j=1}^{d_{k}}\bigg|\frac{f_{j}(z)}{e^{2\pi z_1}}\bigg|^{2}
\end{align*}
%%%%%%%%%%%%%%%%%%%%%%%%%%%%%%%%%%%%%%%%%%%%%%%%%%%%%%%%%%%%%%%%%%%%%%%%%%%%%%%%%%%%%%%%%%%
is a sum of subharmonic functions, and hence, a subharmonic function. So from the maximum principle for subharmonic functions, we deduce that
\begin{align}\label{lem3-eqn1}
\sup_{z\in \overline{U}_{4\pi \slash k(n+1)}(\infty)}\sum_{j=1}^{d_{k}}\bigg|\frac{f_{j}(z)}{e^{2\pi z_1}}\bigg|^{2}=\sup_{z\in \partial U_{4\pi \slash k(n+1)}(\infty)}\sum_{j=1}^{d_{k}}\bigg|\frac{f_{j}(z)}{e^{2\pi z_1}}\bigg|^{2}.
\end{align}
%%%%%%%%%%%%%%%%%%%%%%%%%%%%%%%%%%%%%%%%%%%%%%%%%%%%%%%%%%%%%%%%%%%%%%%%%%%%%%%%%%%%%%%%%%%
We now ascertain the values of $z=(x_1+iy_1,\ldots,x_n+iy_n)\in U_{4\pi\slash k(n+1)}$, for which the term 
\begin{align*}
P(z):=e^{4\pi \mathrm{Re}(z_1)}\big( -2\mathrm{Re}(z_1)-\sum_{i=2}^{n}|z_{i}|^{2}\big)^{k(n+1)}=e^{4\pi x_1}\big( -2x_1-\sum_{i=2}^{n}x_{i}^{2}-\sum_{i=2}^{n}y_{i}^{2}\big)^{k(n+1)}.
\end{align*}
%%%%%%%%%%%%%%%%%%%%%%%%%%%%%%%%%%%%%%%%%%%%%%%%%%%%%%%%%%%%%%%%%%%%%%%%%%%%%%%%%%%%%%%%%%%
attains its maximum.
%%%%%%%%%%%%%%%%%%%%%%%%%%%%%%%%%%%%%%%%%%%%%%%%%%%%%%%%%%%%%%%%%%%%%%%%%%%%%%%%%%%%%%%%%%%

\vspace{0.1cm}
We now compute the derivatives of $P(z)$, with respect to $x_1,\ldots ,x_n$, and $y_1,\ldots,y_n$. 
\begin{align*}
&\frac{\partial P(z)}{\partial x_{1}}=4\pi e^{4\pi x_1}\big( -2x_1-\sum_{i=2}^{n}x_{i}^{2}-\sum_{i=2}^{n}y_{i}^{2}\big)^{k(n+1)}-\\[0.12cm]
&\quad\,\, k(n+1)e^{4\pi x_1}\big( -2x_1-\sum_{i=2}^{n}x_{i}^{2}-\sum_{i=2}^{n}y_{i}^{2}\big)^{k(n+1)-1};\\[0.12cm]
&\frac{\partial P(z)}{\partial x_{i}}=-2k(n+1)x_{i}e^{4\pi x_1}\big( -2x_1-\sum_{i=2}^{n}x_{i}^{2}-\sum_{i=2}^{n}y_{i}^{2}\big)^{k(n+1)-1},\,\,\mathrm{for\,\,all}\,\,2\leq i\leq n;\\[0.12cm]
&\frac{\partial P(z)}{\partial y_{i}}=-2k(n+1)y_{i}e^{4\pi x_1}\big( -2x_1-\sum_{i=2}^{n}x_{i}^{2}-\sum_{i=2}^{n}y_{i}^{2}\big)^{k(n+1)-1},\,\,\mathrm{for\,\,all}\,\,1\leq i\leq n.
\end{align*}
%%%%%%%%%%%%%%%%%%%%%%%%%%%%%%%%%%%%%%%%%%%%%%%%%%%%%%%%%%%%%%%%%%%%%%%%%%%%%%%%%%%%%%%%%%%
Equating the above equations to zero, we deduce that the function $P(z)$ attains a maxima or minima on the set
\begin{align*}
\big\lbrace 2x_1=-k(n+1)\slash 4\pi,z_2=\cdots=z_n=0\big\rbrace=\partial U_{4\pi\slash k(n+1)}(\infty).
\end{align*} 
%%%%%%%%%%%%%%%%%%%%%%%%%%%%%%%%%%%%%%%%%%%%%%%%%%%%%%%%%%%%%%%%%%%%%%%%%%%%%%%%%%%%%%%%%%%
Observing that the function $P(z)$ approaches zero, as $x_{1}$ approaches $-\infty$, we infer that it attains its maxima on the set $\partial U_{4\pi\slash k(n+1)}(\infty)$, which 
together with equation \eqref{lem3-eqn1}, completes the proof of the lemma.
\end{proof}
\end{lem}
%%%%%%%%%%%%%%%%%%%%%%%%%%%%%%%%%%%%%%%%%%%%%%%%%%%%%%%%%%%%%%%%%%%%%%%%%%%%%%%%%%%%%%%%%%%

\vspace{0.1cm}
In the following theorem, we prove Main Theorem \ref{mainthm1}. 
%%%%%%%%%%%%%%%%%%%%%%%%%%%%%%%%%%%%%%%%%%%%%%%%%%%%%%%%%%%%%%%%%%%%%%%%%%%%%%%%%%%%%%%%%%%

\vspace{0.1cm}
\begin{thm}\label{thm4}
With hypothesis as above, for $k\gg 1$, we have the following estimate
\begin{align}\label{thm4:eqn}
\sup_{z\in \overline{X}_{\G}}\big|\bkx(z,z) \big|_{\mathrm{hyp}}=O_{X_{\G}}\big(k^{n+1\slash 2}\big),
\end{align}
%%%%%%%%%%%%%%%%%%%%%%%%%%%%%%%%%%%%%%%%%%%%%%%%%%%%%%%%%%%%%%%%%%%%%%%%%%%%%%%%%%%%%%%%%%%
where the implied constant depends only on $X_{\G}$.
%%%%%%%%%%%%%%%%%%%%%%%%%%%%%%%%%%%%%%%%%%%%%%%%%%%%%%%%%%%%%%%%%%%%%%%%%%%%%%%%%%%%%%%%%%%
\begin{proof}
For a given $\epsilon_{0}>0$, which is independent of $K$, let $\cala_{\epsilon_{0}}:= \calf_{\G_{3}}\backslash U_{\epsilon_{0}}(\infty)$ be a fixed compact subset.  We now estimate
\begin{align*}
&\sup_{z\in\cala_{\epsilon_{0}}}\big|\bkx(\gamma_{31}z,\gamma_{31}z) \big|_{\mathrm{hyp}}.
\end{align*}
%%%%%%%%%%%%%%%%%%%%%%%%%%%%%%%%%%%%%%%%%%%%%%%%%%%%%%%%%%%%%%%%%%%%%%%%%%%%%%%%%%%%%%%%%%%
 
 Adapting arguments from Theorem 4.2 in \cite{anil4}, and applying Proposition \ref{prop2}, we derive
\begin{align} \label{thm4-eqn1}
&\sup_{z\in\cala_{\epsilon_{0}}}\big|\bkx(\gamma_{31}z,\gamma_{31}z) \big|_{\mathrm{hyp}}\leq\ckn+
\frac{\ckn\sin^{2n}(5r_{\epsilon_{0}}\slash 8)}{\sinh^{2n}(r_{\epsilon_0}\slash 4)\cosh^{k(n+1)}(3r_{\epsilon_0}\slash 8)}+\notag\\[0.12cm]
&\frac{\ckn\cosh^{2n}(r_{\epsilon_0}\slash 4)}{(k(n+1)-2n-1)\sinh^{2n}(r_{\epsilon_0}\slash 4)\cosh^{k(n+1)-2n-2}(3r_{\epsilon_{0}}/8)}=O_{X_{\G}}\big(k^{n}\big),\notag\\[0.13cm]&\mathrm{where}\,\,r_{\epsilon_{0}}:=\inf\big\lbrace\dhyp(z,\gamma z)\big|\,z\in \cala_{\epsilon_{0}},\,\,\gamma\in\G_{3}\backslash \lbrace \mathrm{Id}\rbrace \big\rbrace.
\end{align}
%%%%%%%%%%%%%%%%%%%%%%%%%%%%%%%%%%%%%%%%%%%%%%%%%%%%%%%%%%%%%%%%%%%%%%%%%%%%%%%%%%%%%%%%%%%

For $k\gg 1$, it is clear that $U_{4\pi\slash k(n+1)}(\infty)\subset U_{\epsilon_{0}}(\infty)$. So combining equation \eqref{bkiso} with estimate \eqref{bkineq}, and using Lemma \ref{lem3}, we now estimate
\begin{align}\label{thm4-eqn2}
&\sup_{z\in \overline{U}_{4\pi\slash k(n+1)}(\infty)}\big|\bkx(\gamma_{31}z,\gamma_{31}z) \big|_{\mathrm{hyp}}=\notag\\[0.12cm]&\sup_{z\in \partial U_{4\pi\slash k(n+1)}(\infty)}\big| \bkx(\gamma_{31}z,\gamma_{31}z) \big|_{\mathrm{hyp}}\leq\sup_{z\in X_{\G_{3}}}\sum_{\gamma\in\G_3\backslash \G_{3,\infty}}\frac{\ckn}{\cosh^{k(n+1)}(\dhyp(z,\gamma z)\slash 2)}+
\notag\\[0.12cm]&\hspace{4.8cm}\quad\quad\sup_{z\in \partial U_{4\pi\slash k(n+1)}(\infty)}\sum_{\gamma\in\G_{3,\infty}}\frac{\ckn}{\cosh^{k(n+1)}(\dhyp(z,\gamma z)\slash 2)}.
\end{align}
%%%%%%%%%%%%%%%%%%%%%%%%%%%%%%%%%%%%%%%%%%%%%%%%%%%%%%%%%%%%%%%%%%%%%%%%%%%%%%%%%%%%%%%%%%%

We now estimate the two terms on the right hand-side of the above inequality. Adapting arguments from Theorem 4.2 in \cite{anil4}, and applying Proposition \ref{prop2}, we have the following estimate, for the first term on the right hand-side of inequality \eqref{thm4-eqn2}
\begin{align}\label{thm4-eqn3}
&\sup_{z\in X_{\G_{3}}}\sum_{\gamma\in\G_3\backslash \G_{3,\infty}}\frac{\ckn}{\cosh^{k(n+1)}(\dhyp(z,\gamma z)\slash 2)}\leq\notag\\[0.12cm]
&\frac{\ckn\cosh^{2n}(\rx\slash 4)}{(k(n+1)-2n-1)\sinh^{2n}(\rx\slash 4)}+\frac{\ckn\sin^{2n}(5\rx\slash 8)}{\sinh^{2n}(\rx\slash 4)\cosh^{k(n+1)}(3\rx\slash 8)}=O_{X_{\G}}\big(k^{n}\big),
\end{align}
%%%%%%%%%%%%%%%%%%%%%%%%%%%%%%%%%%%%%%%%%%%%%%%%%%%%%%%%%%%%%%%%%%%%%%%%%%%%%%%%%%%%%%%%%%%
where $\rx$ is as defined in equation \eqref{irad}.
%%%%%%%%%%%%%%%%%%%%%%%%%%%%%%%%%%%%%%%%%%%%%%%%%%%%%%%%%%%%%%%%%%%%%%%%%%%%%%%%%%%%%%%%%%%

\vspace{0.1cm}
Combining equations  \eqref{coshreln-1} and \eqref{gamma3infty} with Lemma \ref{lem3}, for any $\gamma\in \G_{3,\infty}$, and $z\in \partial  U_{4\pi\slash k(n+1)}(\infty)$, we infer that 
\begin{align}\label{thm4-eqn4}
&\langle \tilde{z},\tilde{z}\rangle_{H_3}=\langle \tilde{\gamma z},\tilde{\gamma z}\rangle_{H_3}=2\mathrm{Re}(z_1)+\sum_{j=2}^{n}|z_{2}|^{2}=-\frac{k(n+1)}{4\pi};\notag\\[0.12cm]
&\langle\tilde{\gamma z}, \tilde{z}\rangle_{H_3}=2\mathrm{Re}(z_1)+\sum_{j=2}^{n}|z_{2}|^{2}+\frac{-|\tau|^{2}+it}{2}=-\frac{k(n+1)}{4\pi}+\frac{-|\tau|^{2}+it}{2};\notag\\[0.12cm]&\langle\tilde{z}, \tilde{\gamma z}\rangle_{H_3}=2\mathrm{Re}(z_1)+\sum_{j=2}^{n}|z_{2}|^{2}+\frac{-|\tau|^{2}-it}{2}-\frac{k(n+1)}{4\pi}+\frac{-|\tau|^{2}-it}{2},
\end{align}
%%%%%%%%%%%%%%%%%%%%%%%%%%%%%%%%%%%%%%%%%%%%%%%%%%%%%%%%%%%%%%%%%%%%%%%%%%%%%%%%%%%%%%%%%%%
where $\tilde{z}$ and $ \tilde{\gamma z}$ are lifts of $z$ and $\gamma z$ to $\cala^{n+1}$, respectively, and $\tau\in \calo_{K}^{n-1}$, and $t\in\R$.
%%%%%%%%%%%%%%%%%%%%%%%%%%%%%%%%%%%%%%%%%%%%%%%%%%%%%%%%%%%%%%%%%%%%%%%%%%%%%%%%%%%%%%%%%%%

\vspace{0.1cm}
Consider the lattice
\begin{align*}
L:=\bigg\lbrace (\tau,t)\in\C^{n-1}\times \R\bigg|\,\left(\begin{smallmatrix} 1 & -\tau^* U & \frac{-|\tau |^2 + it}{2} \\ 0 & U & \tau \\ 0 & 0 & 1 \end{smallmatrix}\right)\in
\G_{3,\infty},\,U\in\mathrm{SU}(n-1,\calo_{K})\bigg\rbrace\subset \C^{n-1}\times \R.
\end{align*}
%%%%%%%%%%%%%%%%%%%%%%%%%%%%%%%%%%%%%%%%%%%%%%%%%%%%%%%%%%%%%%%%%%%%%%%%%%%%%%%%%%%%%%%%%%%

For a fixed $\tau\in \calo_{K}^{n-1}$, and $t\in\R$, put
\begin{align*}
L_{\tau,t}:=\bigg\lbrace  U\in\mathrm{SU}(n-1,\calo_{K})\bigg|\,\left(\begin{smallmatrix} 1 & -\tau^* U & \frac{-|\tau|^2 + it}{2} \\ 0 & U & \tau \\ 0 & 0 & 1 \end{smallmatrix}\right)\in\G_{3,\infty}\bigg\rbrace.
\end{align*}
%%%%%%%%%%%%%%%%%%%%%%%%%%%%%%%%%%%%%%%%%%%%%%%%%%%%%%%%%%%%%%%%%%%%%%%%%%%%%%%%%%%%%%%%%%%

Since $\mathrm{SU}(n-1,\C)$ is a compact Lie group, we have the following estimate
\begin{align}\label{thm4-eqn5}
\big|L_{\tau,t}\big|=O_{X_{\G}}(1),
\end{align}
%%%%%%%%%%%%%%%%%%%%%%%%%%%%%%%%%%%%%%%%%%%%%%%%%%%%%%%%%%%%%%%%%%%%%%%%%%%%%%%%%%%%%%%%%%%
where the implied constant depends only on $X_{\G}$.
%%%%%%%%%%%%%%%%%%%%%%%%%%%%%%%%%%%%%%%%%%%%%%%%%%%%%%%%%%%%%%%%%%%%%%%%%%%%%%%%%%%%%%%%%%%

\vspace{0.1cm}
We now estimate the second term on the right hand-side of inequality \eqref{thm4-eqn2}. For any $z\in\partial U_{4\pi\slash k(n+1)}(\infty)$, combining the distance formula stated in equation \eqref{coshreln-1} with the relations described in equation \eqref{thm4-eqn4}, and estimate \eqref{thm4-eqn5}, we derive
\begin{align*}
&\sum_{\gamma\in\G_{3,\infty}}\frac{\ckn}{\cosh^{k(n+1)}(\dhyp(z,\gamma z)\slash 2)}\ll_{X_{\G}}\notag\\[0.12cm]&
\sum_{(\tau,t)\in L}\frac{\ckn \big(k(n+1)\slash 4\pi\big)^{k(n+1)}}{\big|k(n+1)\slash 4\pi+|\tau|^{2}\slash 2+it \big|^{k(n+1)}}=\sum_{(\tau,t)\in L}\frac{\ckn \big(k(n+1)\slash 4\pi\big)^{k(n+1)}}{\big((k(n+1)\slash 4\pi+|\tau|^{2}\slash 2)^{2}+t^{2} \big)^{k(n+1)\slash 2}}.
\end{align*}
%%%%%%%%%%%%%%%%%%%%%%%%%%%%%%%%%%%%%%%%%%%%%%%%%%%%%%%%%%%%%%%%%%%%%%%%%%%%%%%%%%%%%%%%%%%

Approximating the summation over the discrete lattice $L$, we derive
\begin{align}\label{thm4-eqn6}
&\sum_{\gamma\in\G_{3,\infty}}\frac{\ckn}{\cosh^{k(n+1)}(\dhyp(z,\gamma z)\slash 2)}\ll_{X_{\G}}\notag\\
&\sum_{u_{n-1}\in\Z}\ldots\sum_{u_1\in\Z}\sum_{v_{n-1}\in\Z}\ldots\sum_{v_1\in\Z}\sum_{t\in \Z}
\frac{\ckn \big(k(n+1)\slash 4\pi\big)^{k(n+1)}}{\big(\big(k(n+1)\slash 4\pi+\sum_{j=1}^{n-1}(u_j^{2}+v_j^2)\big)^{2}+t^{2} \big)^{k(n+1)\slash 2}}.
\end{align}
%%%%%%%%%%%%%%%%%%%%%%%%%%%%%%%%%%%%%%%%%%%%%%%%%%%%%%%%%%%%%%%%%%%%%%%%%%%%%%%%%%%%%%%%%%%

We now estimate the first summation on the right hand-side of inequality \eqref{thm4-eqn6}
\begin{align}\label{thm4-eqn7}
\sum_{t\in \Z}\frac{\ckn \big(k(n+1)\slash 4\pi\big)^{k(n+1)}}{\big(\big(k(n+1)\slash 4\pi+\sum_{j=1}^{n-1}(u_j^{2}+v_j^2)\big)^{2}+t^{2} \big)^{k(n+1)\slash 2}}\leq\notag\\[0.12cm]
\frac{\ckn\big(k(n+1)\slash 4\pi\big)^{k(n+1)}}{\big(k(n+1)\slash 4\pi+\sum_{j=1}^{n-1}(u_j^{2}+v_j^2)\big)^{k(n+1)}}+\notag\\[0.12cm]\int_{0}^{\infty}\frac{2\ckn \big(k(n+1)\slash 4\pi\big)^{k(n+1)}dt}{\big(\big(k(n+1)\slash 4\pi+\sum_{j=1}^{n-1}(u_j^{2}+v_j^2)\big)^{2}+t^{2} \big)^{k(n+1)\slash 2}}.
\end{align}
%%%%%%%%%%%%%%%%%%%%%%%%%%%%%%%%%%%%%%%%%%%%%%%%%%%%%%%%%%%%%%%%%%%%%%%%%%%%%%%%%%%%%%%%%%%

For any $a\in\R_{\geq 1}$, combining formulae 3.251.2 and 8.384.1 from \cite{tables}, we get
\begin{align}\label{formula}
\int_{0}^{\infty}\frac{dx}{\big( 1+x^{2}\big)^{a}}=\frac{\sqrt{\pi}\,\G(a-1\slash 2 )}{\G(a)}.
\end{align}
%%%%%%%%%%%%%%%%%%%%%%%%%%%%%%%%%%%%%%%%%%%%%%%%%%%%%%%%%%%%%%%%%%%%%%%%%%%%%%%%%%%%%%%%%%%

Furthermore, for $a\gg1$, we have the following estimate
\begin{align}\label{formula1}
\int_{0}^{\infty}\frac{dx}{\big( 1+x^{2}\big)^{a}}=\frac{\sqrt{\pi}\,\G(a-1\slash 2 )}{\G(a)}\ll\frac{1}{\sqrt{a}}.
\end{align}
%%%%%%%%%%%%%%%%%%%%%%%%%%%%%%%%%%%%%%%%%%%%%%%%%%%%%%%%%%%%%%%%%%%%%%%%%%%%%%%%%%%%%%%%%%%

Using formula \eqref{formula}, and making the substitution 
\begin{align*}
\theta:=t\slash \big(k(n+1)\slash 4\pi+\sum_{j=1}^{n}(u_j^{2}+v_j^2)\big),
\end{align*}
%%%%%%%%%%%%%%%%%%%%%%%%%%%%%%%%%%%%%%%%%%%%%%%%%%%%%%%%%%%%%%%%%%%%%%%%%%%%%%%%%%%%%%%%%%%
we compute
\begin{align}\label{thm4-eqn8}
&\int_{0}^{\infty}\frac{\ckn \big(k(n+1)\slash 4\pi\big)^{k(n+1)} dt}{\big(\big(k(n+1)\slash 4\pi+\sum_{j=1}^{n-1}(u_j^{2}+v_j^2)\big)^{2}+t^{2} \big)^{k(n+1)\slash 2}}=\notag\\[0.12cm]&\frac{\ckn \big(k(n+1)\slash 4\pi\big)^{k(n+1)} }{\big(k(n+1)\slash 4\pi+\sum_{j=1}^{n-1}(u_j^{2}+v_j^2) \big)^{k(n+1)-1}}\int_{0}^{\infty}\frac{d\theta}{\big(\theta^{2}+1 \big)^{k(n+1)\slash 2}}\ll\notag\\[0.12cm]&
\frac{\ckn \big(k(n+1)\slash 4\pi\big)^{k(n+1)} }{\big(k(n+1)\slash 4\pi+\sum_{j=1}^{n-1}(u_j^{2}+v_j^2) \big)^{k(n+1)-1}}\cdot\frac{\G((k(n+1)-1)\slash 2)}{\G(k(n+1)\slash 2)}.
\end{align}
%%%%%%%%%%%%%%%%%%%%%%%%%%%%%%%%%%%%%%%%%%%%%%%%%%%%%%%%%%%%%%%%%%%%%%%%%%%%%%%%%%%%%%%%%%%

Combining estimates \eqref{thm4-eqn7} and \eqref{thm4-eqn8}, we arrive at the following estimate
\begin{align}\label{thm4-eqn9}
&\sum_{u_{n-1}\in\Z}\ldots\sum_{u_1\in\Z}\sum_{v_{n-1}\in\Z}\ldots\sum_{v_1\in\Z}\sum_{t\in \Z}\frac{\ckn \big(k(n+1)\slash 4\pi\big)^{k(n+1)}}{\big(\big(k(n+1)\slash 4\pi+\sum_{j=1}^{n-1}(u_j^{2}+v_j^2)\big)^{2}+t^{2} \big)^{k(n+1)\slash 2}}\ll\notag\\[0.12cm]&\sum_{u_{n-1}\in\Z}\ldots\sum_{u_1\in\Z}\sum_{v_{n-1}\in\Z}\ldots\sum_{v_1\in\Z}\frac{\ckn\big(k(n+1)\slash 4\pi\big)^{k(n+1)}}{\big(k(n+1)\slash 4\pi+\sum_{j=1}^{n-1}(u_j^{2}+v_j^2)\big)^{k(n+1)}}+\notag\\[0.12cm]&\sum_{u_{n-1}\in\Z}\ldots\sum_{u_1\in\Z}\sum_{v_{n-1}\in\Z}\ldots\sum_{v_1\in\Z}\frac{\ckn \big(k(n+1)\slash 4\pi\big)^{k(n+1)} }{\big(k(n+1)\slash 4\pi+\sum_{j=1}^{n-1}(u_j^{2}+v_j^2) \big)^{k(n+1)-1}}\cdot\frac{\G((k(n+1)-1)\slash 2)}{\G(k(n+1)\slash 2)}.
\end{align}
%%%%%%%%%%%%%%%%%%%%%%%%%%%%%%%%%%%%%%%%%%%%%%%%%%%%%%%%%%%%%%%%%%%%%%%%%%%%%%%%%%%%%%%%%%%

We now estimate the second summation in the variable $v_1$. From similar arguments as the ones used to prove estimates \eqref{thm4-eqn7} and \eqref{thm4-eqn8}, we derive 
\begin{align}\label{thm4-eqn10}
\sum_{v_1\in\Z}\frac{\ckn\big(k(n+1)\slash 4\pi\big)^{k(n+1)}}{\big(k(n+1)\slash 4\pi+\sum_{j=1}^{n-1}(u_j^{2}+v_j^2)\big)^{k(n+1)}}\ll
\frac{\ckn\big(k(n+1)\slash 4\pi\big)^{k(n+1)}}{\big(k(n+1)\slash 4\pi+\sum_{j=2}^{n-1}(u_j^{2}+v_j^2)+u_{1}^{2}\big)^{k(n+1)}}+\notag\\[0.12cm]\frac{\ckn\big(k(n+1)\slash 4\pi\big)^{k(n+1)}}{\big(k(n+1)\slash 4\pi+\sum_{j=2}^{n-1}(u_j^{2}+v_j^2)+u_{1}^{2}\big)^{k(n+1)-1\slash 2}}\cdot\frac{\G(k(n+1)-1)}{\G(k(n+1)-1\slash 2)},
\end{align}
%%%%%%%%%%%%%%%%%%%%%%%%%%%%%%%%%%%%%%%%%%%%%%%%%%%%%%%%%%%%%%%%%%%%%%%%%%%%%%%%%%%%%%%%%%%
and
\begin{align}\label{thm4-eqn11}
\sum_{v_1\in\Z}\frac{\ckn \big(k(n+1)\slash 4\pi\big)^{k(n+1)} }{\big(k(n+1)\slash 4\pi+\sum_{j=1}^{n-1}(u_j^{2}+v_j^2) \big)^{k(n+1)-1}}\ll \notag\\[0.12cm]
\frac{\ckn\big(k(n+1)\slash 4\pi\big)^{k(n+1)}}{\big(k(n+1)\slash 4\pi+\sum_{j=2}^{n-1}(u_j^{2}+v_j^2)+u_{1}^{2}\big)^{k(n+1)-1}} +\notag\\[0.12cm]
\frac{\ckn \big(k(n+1)\slash 4\pi\big)^{k(n+1)} }{\big(k(n+1)\slash 4\pi+\sum_{j=2}^{n-1}(u_j^{2}+v_j^2)+u_{1}^{2} \big)^{k(n+1)-3\slash 2}}\cdot\frac{\G(k(n+1)-3\slash 2)}{\G(k(n+1)-1)}.
\end{align}
%%%%%%%%%%%%%%%%%%%%%%%%%%%%%%%%%%%%%%%%%%%%%%%%%%%%%%%%%%%%%%%%%%%%%%%%%%%%%%%%%%%%%%%%%%%

Combining estimate \eqref{thm4-eqn6}, with estimates \eqref{thm4-eqn8}--\eqref{thm4-eqn11}, and applying estimate \eqref{formula1}, we arrive at the following estimate
\begin{align}\label{thm4-eqn12}
&\sum_{\gamma\in\G_{3,\infty}}\frac{\ckn}{\cosh^{k(n+1)}(\dhyp(z,\gamma z)\slash 2)}\ll_{X_{\G}}\notag\\
& \sum_{u_{n-1}\in\Z}\ldots\sum_{u_1\in\Z}\sum_{v_{n-1}\in\Z}\ldots\sum_{v_2\in\Z}\frac{\ckn \big(k(n+1)\slash 4\pi\big)^{k(n+1)} }{\big(k(n+1)\slash 4\pi+\sum_{j=2}^{n-1}(u_j^{2}+v_j^2)+u_{1}^{2} \big)^{k(n+1)}}+
\notag\\
&\sum_{u_{n-1}\in\Z}\ldots\sum_{u_1\in\Z}\sum_{v_{n-1}\in\Z}\ldots\sum_{v_2\in\Z}\Bigg( \sum_{\alpha=0}^{2}\frac{\ckn \big(k(n+1)\slash 4\pi\big)^{k(n+1)-\alpha\slash 2} }{\big(k(n+1)\slash 4\pi+\sum_{j=2}^{n-1}(u_j^{2}+v_j^2)+u_{1}^{2} \big)^{k(n+1)-(\alpha+1)\slash 2}}\Bigg) .
\end{align}
%%%%%%%%%%%%%%%%%%%%%%%%%%%%%%%%%%%%%%%%%%%%%%%%%%%%%%%%%%%%%%%%%%%%%%%%%%%%%%%%%%%%%%%%%%%

Repeating the same process for computing the sums in the the other variables, and applying Proposition \ref{prop2}, we arrive at the following estimate
\begin{align}\label{thm4-eqn13}
\sup_{z\in\partial U_{4\pi\slash k(n+1)}(\infty)}\sum_{\gamma\in\G_{3,\infty}}\frac{\ckn}{\cosh^{k(n+1)}(\dhyp(z,\gamma z)\slash 2)}\ll_{X_{\G}}\ckn\sqrt{k}=O_{X_{\G}}(k^{n+1\slash 2}). 
\end{align}
%%%%%%%%%%%%%%%%%%%%%%%%%%%%%%%%%%%%%%%%%%%%%%%%%%%%%%%%%%%%%%%%%%%%%%%%%%%%%%%%%%%%%%%%%%%

Combining estimates \eqref{thm4-eqn2}, \eqref{thm4-eqn3}, and \eqref{thm4-eqn13}, we arrive at the following estimate
\begin{align}\label{thm4-eqn14}
&\sup_{z\in \overline{U}_{4\pi\slash k(n+1)}(\infty)}\big|\bkx(\gamma_{31}z,\gamma_{31}z) \big|_{\mathrm{hyp}}=O_{X_{\G}}(k^{n+1\slash 2}).
\end{align} 
%%%%%%%%%%%%%%%%%%%%%%%%%%%%%%%%%%%%%%%%%%%%%%%%%%%%%%%%%%%%%%%%%%%%%%%%%%%%%%%%%%%%%%%%%%%

We now estimate
\begin{align*}
\sup_{z\in U_{\epsilon_{0}}(\infty)\backslash U_{4\pi\slash k(n+1)}(\infty)}\big|\bkx(\gamma_{31}z,\gamma_{31}z) \big|_{\mathrm{hyp}}\leq
\sup_{z\in X_{\G_{3}}}\sum_{\gamma\in\G_3\backslash \G_{3,\infty}}\frac{\ckn}{\cosh^{k(n+1)}(\dhyp(z,\gamma z)\slash 2)}+\notag\\[0.12cm]\sup_{z\in U_{\epsilon_{0}}(\infty)\backslash U_{4\pi\slash k(n+1)}(\infty)}\sum_{\gamma\in\G_{3,\infty}}\frac{\ckn}{\cosh^{k(n+1)}(\dhyp(z,\gamma z)\slash 2)}.
\end{align*}
%%%%%%%%%%%%%%%%%%%%%%%%%%%%%%%%%%%%%%%%%%%%%%%%%%%%%%%%%%%%%%%%%%%%%%%%%%%%%%%%%%%%%%%%%%%

From relations described in equation \eqref{thm4-eqn4}, for any $\gamma\in \G_{3,\infty}$, and $z\in U_{\epsilon_{0}}(\infty)\backslash U_{2\pi\slash k}(\infty)$, we find that
\begin{align*}
&\langle \tilde{z},\tilde{z}\rangle_{H_3}=\langle \tilde{\gamma z},\tilde{\gamma z}\rangle_{H_3}=2\mathrm{Re}(z_1)+\sum_{j=2}^{n}|z_{2}|^{2}=2x_1;\notag\\
&\langle\tilde{\gamma z}, \tilde{z}\rangle_{H_3}=2x_1+\frac{-|\tau|^{2}+it}{2};\,\,\langle\tilde{z}, \tilde{\gamma z}\rangle_{H_3}=2x_1+\frac{-|\tau|^{2}-it}{2},
\end{align*}
%%%%%%%%%%%%%%%%%%%%%%%%%%%%%%%%%%%%%%%%%%%%%%%%%%%%%%%%%%%%%%%%%%%%%%%%%%%%%%%%%%%%%%%%%%%
where $\tilde{z}$ and $ \tilde{\gamma z}$ are lifts of $z$ and $\gamma z$ to $\cala^{n+1}$, respectively, and $\tau\in \calo_{K}^{n-1}$, and $t\in\R$.
%%%%%%%%%%%%%%%%%%%%%%%%%%%%%%%%%%%%%%%%%%%%%%%%%%%%%%%%%%%%%%%%%%%%%%%%%%%%%%%%%%%%%%%%%%%

From the arguments employed to prove estimate \eqref{thm4-eqn12}, we deduce that 
\begin{align}\label{thm4-eqn15}
&\sup_{z\in U_{\epsilon_{0}}(\infty)\backslash U_{4\pi\slash k(n+1)}(\infty)}\sum_{\gamma\in\G_{3,\infty}}\frac{\ckn}{\cosh^{k(n+1)}(\dhyp(z,\gamma z)\slash 2)}\ll_{X_{\G}}\notag\\
&\sum_{u_{n-1}\in\Z}\ldots\sum_{u_1\in\Z}\sum_{v_{n-1}\in\Z}  \ldots\sum_{v_2\in\Z}\Bigg( \sum_{\alpha=0}^{2}\frac{\ckn \big(2x_1\big)^{k(n+1)-\alpha\slash 2} }{\big(2x_1+\sum_{j=2}^{n-1}(u_j^{2}+v_j^2)+u_{1}^{2} \big)^{k(n+1)-(\alpha+1)\slash 2}}\Bigg)=\notag\\[0.125cm]&\hspace{8.2cm} O_{X_{\G}}\big(\ckn\sqrt{x_1}\big)=O_{X_{\G}}\big(k^{n+1\slash 2}\big). 
\end{align}
%%%%%%%%%%%%%%%%%%%%%%%%%%%%%%%%%%%%%%%%%%%%%%%%%%%%%%%%%%%%%%%%%%%%%%%%%%%%%%%%%%%%%%%%%%%

Combining estimates \eqref{thm4-eqn1}, \eqref{thm4-eqn14}, and \eqref{thm4-eqn15}, completes the proof of the theorem.
\end{proof}
\end{thm}
%%%%%%%%%%%%%%%%%%%%%%%%%%%%%%%%%%%%%%%%%%%%%%%%%%%%%%%%%%%%%%%%%%%%%%%%%%%%%%%%%%%%%%%%%%%
%%%%%%%%%%%%%%%%%%%%%%%%%%%%%%%%%%%%%%%%%%%%%%%%%%%%%%%%%%%%%%%%%%%%%%%%%%%%%%%%%%%%%%%%%%%

\vspace{0.2cm}
\subsection{Estimates of the Bergman metric}\label{subsec-3.2}
%%%%%%%%%%%%%%%%%%%%%%%%%%%%%%%%%%%%%%%%%%%%%%%%%%%%%%%%%%%%%%%%%%%%%%%%%%%%%%%%%%%%%%%%%%%
We now prove Main Theorem \ref{mainthm2}, and we continue to work with a fixed fundamental domain $(1,0,\ldots,0)^{t}\in\calf_{\G}$.
%%%%%%%%%%%%%%%%%%%%%%%%%%%%%%%%%%%%%%%%%%%%%%%%%%%%%%%%%%%%%%%%%%%%%%%%%%%%%%%%%%%%%%%%%%%

\vspace{0.1cm}
Furthermore, for brevity of notation for ensuing computations, we denote $\bkn(z,z)$ by $\bkn(z)$. For $z\in\Bn$, from equation \eqref{bkseries}, we infer that 
\begin{align}\label{bkseries-1}
&\bkn(z)=\sum_{\gamma\in \G}\frac{\ckn}{\big(\big(\sum_{j=1}^{n}\overline{c_{j}z_{j}}\big)+\overline{D}-\big(\sum_{i=1}^{n}\big(\sum_{j=1}^{n}\overline{a_{ij}z_{j}}\big)+\overline{b_{i}}\big)z_{i}\big)^{k(n+1)}},\notag\\[0.13cm] 
&\mathrm{where}\,\gamma=\left(\begin{smallmatrix}A &B\\C&D  \end{smallmatrix}\right)\in\G\,\,\mathrm{and}\,\, A:=(a_{ij})_{1\leq i,j\leq n}\in M_{n\times n}(\C),\notag\\[0.12cm] &B:=(b_1,\ldots,b_n)^{t}\in M_{n\times 1}(\C),\, C:=(c_1,\ldots,c_n)\in M_{1\times n}(\C), \mathrm{and}\,\,D\in \C.
\end{align}
%%%%%%%%%%%%%%%%%%%%%%%%%%%%%%%%%%%%%%%%%%%%%%%%%%%%%%%%%%%%%%%%%%%%%%%%%%%%%%%%%%%%%%%%%%%

\vspace{0.1cm}
\begin{lem}\label{lem5}
With notation as above, for any $z\in \calf_{\G}$, we have
\begin{align*}
\berkvol(z)=\bigg(-\frac{i}{2\pi}\bigg)^{n}\mathrm{det}(M(z))dz_1\wedge d\overline{z}_{1}\wedge\ldots\wedge dz_n\wedge d\overline{z}_n
\end{align*}
%%%%%%%%%%%%%%%%%%%%%%%%%%%%%%%%%%%%%%%%%%%%%%%%%%%%%%%%%%%%%%%%%%%%%%%%%%%%%%%%%%%%%%%%%%%
where $M(z):=\big( m_{ij}(z)\big)_{1\leq i,j\leq n }$, and for $1\leq i\not=j\leq n$, we have
\begin{align}\label{lem5:eqn1}
m_{ij}(z):=-\frac{k(n+1) z_{j}\overline{z}_{i}}{\big( 1-|z|^{2}\big)^{2}}+\frac{1}{\bkn(z)}\frac{\partial^{2}\bkn(z)}{\partial z_{i}\partial\overline{z}_{j}}-
\notag\\[0.12cm]\frac{1}{\big(\bkn(z)\big)^{2}}\frac{\partial\bkn(z)}{\partial z_{i}}\frac{\partial\bkn(z)}{\partial\overline{z}_{j}},
\end{align}
%%%%%%%%%%%%%%%%%%%%%%%%%%%%%%%%%%%%%%%%%%%%%%%%%%%%%%%%%%%%%%%%%%%%%%%%%%%%%%%%%%%%%%%%%%%
and 
\begin{align}\label{lem5:eqn2}
m_{ii}(z):=-\frac{k(n+1)\big(1-\sum_{\substack{j=1\\j\not=i}}^{n} |z_{j}|^{2} \big)}{\big( 1-|z|^{2}\big)^{2}}+\frac{1}{\bkn(z)}\frac{\partial^{2}\bkn(z)}{\partial z_{i}\partial\overline{z}_{i}}-\notag\\[0.12cm]\frac{1}{\big(\bkn(z)\big)^{2}}\bigg|\frac{\partial\bkn(z)}{\partial z_{i}}\bigg|^{2}.
\end{align}
%%%%%%%%%%%%%%%%%%%%%%%%%%%%%%%%%%%%%%%%%%%%%%%%%%%%%%%%%%%%%%%%%%%%%%%%%%%%%%%%%%%%%%%%%%%
\begin{proof}
The lemma follows directly from equation \eqref{bergman-metric-1}.
\end{proof}
\end{lem}
%%%%%%%%%%%%%%%%%%%%%%%%%%%%%%%%%%%%%%%%%%%%%%%%%%%%%%%%%%%%%%%%%%%%%%%%%%%%%%%%%%%%%%%%%%%

\vspace{0.1cm}
In the following proposition, we derive a lower bound for the Bergman kernel.
%%%%%%%%%%%%%%%%%%%%%%%%%%%%%%%%%%%%%%%%%%%%%%%%%%%%%%%%%%%%%%%%%%%%%%%%%%%%%%%%%%%%%%%%%%%

\vspace{0.1cm}
\begin{prop}\label{prop9}
With notation as above, for $k\gg1$ and $z=(z_1,\ldots, z_{n})^{t}\in U_{2\pi/ (k(n+1))}(\infty)$, we have the following lower bound
\begin{align}\label{prop9:eqn1}
\big|\bkn(z)\big|_{\mathrm{pet}}\gg_{X_{\G}}\frac{k^{n}}{\sqrt{-2\mathrm{Re}(z_1)-\sum_{i=2}^{n}|z_{i}|^{2}}};
\end{align}
%%%%%%%%%%%%%%%%%%%%%%%%%%%%%%%%%%%%%%%%%%%%%%%%%%%%%%%%%%%%%%%%%%%%%%%%%%%%%%%%%%%%%%%%%%%
when $z=(z_1,\ldots, z_{n})^{t}\in \mathcal{F}_{\G_{3}}\backslash U_{2\pi/ (k(n+1))}(\infty)$, we have the following lower bound
\begin{align}\label{prop9:eqn2}
\big|\bkn(z)\big|_{\mathrm{pet}}\gg_{X_{\G}}\bigg(\frac{ -2\mathrm{Re}(z_1)-\sum_{i=2}^{n}|z_{i}|^{2}}{k}\bigg)^{(n-1)/2} k^{n-1/2},
\end{align}
%%%%%%%%%%%%%%%%%%%%%%%%%%%%%%%%%%%%%%%%%%%%%%%%%%%%%%%%%%%%%%%%%%%%%%%%%%%%%%%%%%%%%%%%%%%
and the implied constants in both estimates  depends only on $X_{\G}$.
%%%%%%%%%%%%%%%%%%%%%%%%%%%%%%%%%%%%%%%%%%%%%%%%%%%%%%%%%%%%%%%%%%%%%%%%%%%%%%%%%%%%%%%%%%%
\begin{proof}
As in the proof of Theorem \ref{thm4}, for any $z\in  \mathcal{F}_{\G_{3}}$, we have the following estimate
\begin{align}\label{prop9-eqn0}
&\sum_{\gamma\in\G_3\backslash \Gamma_{3,\infty}}\frac{\ckn}{\cosh^{k(n+1)}(\dhyp(z,\gamma z)\slash 2)}\ll_{X_{\G}}\frac{\ckn\sin^{2n}(5r_{X_{\Gamma}}\slash 8)}{\sinh^{2n}(r_{X_{\Gamma}}\slash 4)\cosh^{k(n+1)}(3r_{X_{\Gamma}}\slash 8)}+\notag\\[0.12cm] & \frac{\ckn\cosh^{2n}(r_{X_{\Gamma}}\slash 4)}{(k(n+1)-2n-1)\sinh^{2n}(r_{X_{\Gamma}}\slash 4)\cosh^{k(n+1)-2n-2}(3r_{X_{\Gamma}}/8)}\ll_{X_{\G}}\frac{k^{n}}{\cosh^{k(n+1)}(3r_{X_{\Gamma}}/8)},
\notag\\[0.13cm]&\mathrm{where}\,\,r_{X_{\Gamma}}:=\inf\big\lbrace \mathrm{d}_{\mathrm{hyp}}(z,\gamma z)\big|\,z\in  \mathcal{F}_{\G_{3}},\,\gamma\in\Gamma_{3}\backslash \Gamma_{3,\infty} \big\rbrace.
\end{align}
%%%%%%%%%%%%%%%%%%%%%%%%%%%%%%%%%%%%%%%%%%%%%%%%%%%%%%%%%%%%%%%%%%%%%%%%%%%%%%%%%%%%%%%%%%%
With notation as in Theorem \ref{thm4}, and using equation \eqref{thm4-eqn4}, we have
\begin{align}\label{prop9-eqn1}
&\sum_{\gamma\in \Gamma_{3,\infty}}\frac{\ckn \langle \tilde{z},\tilde{ z}\rangle_{H_{3}}^{k(n+1)}}{\langle \tilde{z},\tilde{\gamma z}\rangle_{H_{3}}^{k(n+1)}}=\sum_{(\tau,L)\in L}
\frac{\ckn \big(-2\mathrm{Re}(z_1)-\sum_{i=2}^{n}|z_{i}|^{2} \big)^{k(n+1)}}{\big(-2\mathrm{Re}(z_1)-\sum_{i=2}^{n}|z_{i}|^{2} +\frac{|\tau|^{2}}{2}+i\frac{t}{2}\big)^{k(n+1)}}.
\end{align}
%%%%%%%%%%%%%%%%%%%%%%%%%%%%%%%%%%%%%%%%%%%%%%%%%%%%%%%%%%%%%%%%%%%%%%%%%%%%%%%%%%%%%%%%%%%
Put 
$$
\alpha(z):=-2\mathrm{Re}(z_1)-\sum_{i=2}^{n}|z_{i}|^{2}. 
$$
%%%%%%%%%%%%%%%%%%%%%%%%%%%%%%%%%%%%%%%%%%%%%%%%%%%%%%%%%%%%%%%%%%%%%%%%%%%%%%%%%%%%%%%%%%%
From Remark \ref{rem-2.1}, and using Poisson summation formula, we have
\begin{align}\label{prop9-eqn2}
\sum_{(\tau,L)\in L}
\frac{\big(-2\mathrm{Re}(z_1)-\sum_{i=2}^{n}|z_{i}|^{2} \big)^{k(n+1)}}{\big(-2\mathrm{Re}(z_1)-\sum_{i=2}^{n}|z_{i}|^{2} +\frac{|\tau|^{2}}{2}+i\frac{t}{2}\big)^{k(n+1)}}\gg_{X_{\Gamma}}
\sum_{\tau\in\mathbb{Z}^{2n-2}}\sum_{t\in\mathbb{Z}}\frac{\alpha(z)^{k(n+1)}}{\big(\alpha(z)+|\tau|^{2}+i t\big)^{k(n+1)}}=\notag\\[0.12cm]\sum_{\tau\in\mathbb{Z}^{2n-2}}\sum_{m\in\mathbb{Z}_{\geq 1}}\frac{2\pi \alpha(z)^{k(n+1)}e^{-2\pi m(|\tau|^{2}+\alpha(z))}(2\pi m)^{k(n+1)}}{(k(n+1)-1)!}.
\end{align}
%%%%%%%%%%%%%%%%%%%%%%%%%%%%%%%%%%%%%%%%%%%%%%%%%%%%%%%%%%%%%%%%%%%%%%%%%%%%%%%%%%%%%%%%%%%
Using integral test for the lattice $\tau\in\mathbb{Z}^{2n-2}$, we derive
\begin{align}\label{prop9-eqn3}
&\sum_{\tau\in\mathbb{Z}^{2n-2}}\sum_{m\in\mathbb{Z}_{\geq 1}}\frac{2\pi \alpha(z)^{k(n+1)}e^{-2\pi m (|\tau|^{2}+\alpha(z))}(2\pi m)^{k(n+1)}}{(k(n+1)-1)!}\gg_{X_{\Gamma}}\notag\\[0.12cm]&\sum_{m\in\mathbb{Z}_{\geq 1}}\bigg(\int_{0}^{\infty}e^{-2\pi mx^{2}}dx\bigg)^{2n-2}\cdot\frac{2\pi \alpha(z)^{k(n+1)}e^{-2\pi m \alpha(z)}(2\pi m)^{k(n+1)}}{(k(n+1)-1)!}\gg_{X_{\Gamma}}\notag\\[0.12cm]&\sum_{m\in\mathbb{Z}_{\geq 1}}\frac{(2\pi \alpha (z))^{k(n+1)}m^{k(n+1)-(n+1)/2}e^{-2\pi m \alpha(z)}}{(k(n+1)-1)!}.
\end{align}
%%%%%%%%%%%%%%%%%%%%%%%%%%%%%%%%%%%%%%%%%%%%%%%%%%%%%%%%%%%%%%%%%%%%%%%%%%%%%%%%%%%%%%%%%%%
Using integral test for the summation on $m\in\mathbb{Z}_{\geq }1$, we now compute
\begin{align}\label{prop9-eqn4}
&\sum_{m\in\mathbb{Z}_{\geq 1}}\frac{(2\pi \alpha (z))^{k(n+1)}m^{k(n+1)-(n+1)/2}e^{-2\pi m \alpha(z)}}{(k(n+1)-1)!}\geq\notag\\[0.12cm]&\frac{(2\pi \alpha (z))^{(n+1)/2}}{(k(n+1)-1)!}
\int_{1}^{\infty}(2\pi m \alpha (z))^{k(n+1)-(n+1)/2}e^{-2\pi \alpha(z)}dm\geq \notag\\[0.12cm]&\frac{(2\pi \alpha (z))^{(n-1)/2}}{(k(n+1)-1)!}\int_{2\pi \alpha (z)}^{\infty}\theta(z)^{k(n+1)-(n+1)/2}e^{-\theta}d\theta.
\end{align}
%%%%%%%%%%%%%%%%%%%%%%%%%%%%%%%%%%%%%%%%%%%%%%%%%%%%%%%%%%%%%%%%%%%%%%%%%%%%%%%%%%%%%%%%%%%
For $z=(z_1,\ldots, z_{n})^{t}\in U_{2\pi/ (k(n+1))}(\infty)$, since $2\pi \alpha(z)>k(n+1)$, combining estimates \eqref{prop9-eqn1}--\eqref{prop9-eqn4} with Proposition \ref{prop2}, using Sterling's approximation, we have the following lower bound
\begin{align}\label{prop9-eqn5}
\sum_{\gamma\in \Gamma_{3,\infty}}\frac{\ckn \langle \tilde{z},\tilde{ z}\rangle_{H_{3}}^{k(n+1)}}{\langle \tilde{z},\tilde{\gamma z}\rangle_{H_{3}}^{k(n+1)}}
\gg_{X_{\Gamma}}\frac{k^{n}}{\sqrt{-2\mathrm{Re}(z_1)-\sum_{i=2}^{n}|z_{i}|^{2}}}.
\end{align}
%%%%%%%%%%%%%%%%%%%%%%%%%%%%%%%%%%%%%%%%%%%%%%%%%%%%%%%%%%%%%%%%%%%%%%%%%%%%%%%%%%%%%%%%%%%
Similarly, for $z=(z_1,\ldots, z_{n})^{t}\in \mathcal{F}_{\Gamma_{3}}\backslash U_{2\pi/ (k(n+1))}(\infty)$, since $k(n+1)<2\pi \alpha(z)$, combining estimates \eqref{prop9-eqn1}--\eqref{prop9-eqn4} with Proposition \ref{prop2}, using Sterling's approximation, we have the following lower bound
\begin{align}\label{prop9-eqn6}
\sum_{\gamma\in \Gamma_{3,\infty}}\frac{\ckn \langle \tilde{z},\tilde{ z}\rangle_{H_{3}}^{k(n+1)}}{\langle \tilde{z},\tilde{\gamma z}\rangle_{H_{3}}^{k(n+1)}}
\gg_{X_{\Gamma}} \bigg(\frac{ -2\mathrm{Re}(z_1)-\sum_{i=2}^{n}|z_{i}|^{2}}{k}\bigg)^{(n-1)/2} k^{n-1/2},
\end{align}
%%%%%%%%%%%%%%%%%%%%%%%%%%%%%%%%%%%%%%%%%%%%%%%%%%%%%%%%%%%%%%%%%%%%%%%%%%%%%%%%%%%%%%%%%%%
For $k\gg1$, combining estimate \eqref{prop9-eqn0} with estimates \eqref{prop9-eqn5} and \eqref{prop9-eqn6}, completes the proof of the proposition. 
\end{proof}
\end{prop}
%%%%%%%%%%%%%%%%%%%%%%%%%%%%%%%%%%%%%%%%%%%%%%%%%%%%%%%%%%%%%%%%%%%%%%%%%%%%%%%%%%%%%%%%%%%

\vspace{0.1cm}
We now estimate the derivatives of the Bergman kernel, using which we estimate the determinant of the matrix $M(z)$. The following proposition, is an extension of Proposition 3.8 
from \cite{anil7}.
%%%%%%%%%%%%%%%%%%%%%%%%%%%%%%%%%%%%%%%%%%%%%%%%%%%%%%%%%%%%%%%%%%%%%%%%%%%%%%%%%%%%%%%%%%%

\vspace{0.1cm}
\begin{prop}\label{prop8}
With notation as above, for $k\gg 1$, and $1\leq i,j \leq n$, and a fixed $a\in \mathbb{R}$ with $-2n\leq a\leq 2n$, we have the following estimate
\begin{align}\label{prop7:eqn}
\sup_{z\in \overline{\calf}_{\G}}\big(1-|z|^{2}\big)^{k(n+1)-a}\bigg|\frac{\partial^{2} \bkn(z)}{\partial z_{i}\partial \overline{z}_{j}} \bigg|=O_{X_{\G}}\big(k^{n+a+9/2} \big),
\end{align}
%%%%%%%%%%%%%%%%%%%%%%%%%%%%%%%%%%%%%%%%%%%%%%%%%%%%%%%%%%%%%%%%%%%%%%%%%%%%%%%%%%%%%%%%%%%
where the implied constant depends only on $X_{\G}$.
%%%%%%%%%%%%%%%%%%%%%%%%%%%%%%%%%%%%%%%%%%%%%%%%%%%%%%%%%%%%%%%%%%%%%%%%%%%%%%%%%%%%%%%%%%%%
\begin{proof}
Let $\lbrace f_1,\ldots, f_{d_{k}}\rbrace$ denote an orthonormal basis of $\skn$. For any $z\in \calf_{\G}$, and $1\leq i,j \leq n$, and a fixed $a\in \mathbb{R}$ with $-2n\leq a\leq 2n$,  from the definition of the Bergman kernel, and the Cauchy-Schwartz inequality, we have 
\begin{align}\label{prop7-eqn0}
\Bigg|\frac{\partial^{2} \bkn(z)}{\partial z_{i}\partial \overline{z}_{j}}\Bigg| =\Bigg|\sum_{l=1}^{d_{k}}\frac{\partial f_{l}(z)}{\partial z_{i}}\frac{\partial \overline{f_{l}(z)}}{\partial \overline{z}_{j}}\Bigg|\leq \sqrt{\sum_{l=1}^{d_{k}}\bigg|\frac{\partial f_{l}(z)}{\partial z_{i}}\bigg|^{2}}\cdot \sqrt{\sum_{l=1}^{d_{k}}\bigg|\frac{\partial f_{l}(z)}{\partial z_{j}}\bigg|^{2}}\leq\notag\\[0.12cm]
\sup_{\substack{z\in \overline{\calf}_{\G},\\i\leq i\leq n}}\sum_{l=1}^{d_{k}}\bigg|\frac{\partial f_{l}(z)}{\partial z_{i}}\bigg|^{2}=
\sup_{\substack{z\in \overline{\calf}_{\G},\\i\leq i\leq n}}\Bigg|\frac{\partial^{2} \bkn(z)}{\partial z_{i}\partial \overline{z}_{i}}\Bigg| .
\end{align}
%%%%%%%%%%%%%%%%%%%%%%%%%%%%%%%%%%%%%%%%%%%%%%%%%%%%%%%%%%%%%%%%%%%%%%%%%%%%%%%%%%%%%%%%%%%%
For any $1\leq i\leq n$, the following series 
\begin{align*}
\Bigg|\frac{\partial^{2} \bkn(z)}{\partial z_{i}\partial \overline{z}_{i}}\Bigg|=\sum_{l=1}^{d_{k}}\bigg|\frac{\partial f_{l}(z)}{\partial z_{i}}\bigg|^{2}
\end{align*}
%%%%%%%%%%%%%%%%%%%%%%%%%%%%%%%%%%%%%%%%%%%%%%%%%%%%%%%%%%%%%%%%%%%%%%%%%%%%%%%%%%%%%%%%%%%%
is a subharmonic function. 
%%%%%%%%%%%%%%%%%%%%%%%%%%%%%%%%%%%%%%%%%%%%%%%%%%%%%%%%%%%%%%%%%%%%%%%%%%%%%%%%%%%%%%%%%%%%

\vspace{0.1cm}
From brevity of notation, set $U_{k,a}:=U_{4\pi/(k(n+1)-a)}(\infty)$, and $\overline{U}_{k,a}:=\overline{U}_{4\pi/(k(n+1)-a)}(\infty)$. Hence, adapting arguments from Lemma \ref{lem3} and Theorem \ref{thm4}, we deduce that 
\begin{align*}
\sup_{z\in \overline{\calf}_{\G}}\big(1-|z|^{2}\big)^{k(n+1)-a}\Bigg|\frac{\partial^{2} \bkn(z)}{\partial z_{i}\partial \overline{z}_{i}}\Bigg|=\sup_{z\in \gamma_{31}(\partial U_{k,a})}\frac{\big(1-|z|^{2}\big)^{k(n+1)+2}}{\big(1-|z|^{2}\big)^{a+2}}\Bigg|\frac{\partial^{2} \bkn(z)}{\partial z_{i}\partial \overline{z}_{i}}\Bigg|.
\end{align*}
%%%%%%%%%%%%%%%%%%%%%%%%%%%%%%%%%%%%%%%%%%%%%%%%%%%%%%%%%%%%%%%%%%%%%%%%%%%%%%%%%%%%%%%%%%%%
For any $1\leq i \leq n$, $\gamma\in \G$, and $z\in \gamma_{31}(\partial U_{k,a})$, we now compute
\begin{align}\label{prop7-eqn1}
&\frac{\partial^{2} \bkn(z)}{\partial z_{i}\partial \overline{z}_{i}} =\sum_{\gamma\in \G}\frac{\ckn k(n+1)a_{ii}}{\big(Cz+D-\sum_{l=1}^{n}\big(\big(\sum_{m=1}^{n}a_{lm}z_{m}\big)+b_{l}\big)\overline{z_{l}}\big)^{k(n+1)+1}}-\notag\\[0.12cm]&\sum_{\gamma\in \G}\frac{\ckn k(n+1)(k(n+1)+1)\big(\big(\sum_{l=1}^{n}a_{il}z_{l}\big) +b_{i}\big)\big( c_{i}-\big(\sum_{l=1}^{n}a_{li}\overline{z_{li}} \big)\big)}{\big(Cz+D-\sum_{l=1}^{n}\big(\big(\sum_{m=1}^{n}a_{lm}z_{m}\big)+b_{l}\big)\overline{z_{l}}\big)^{k(n+1)+2}}.
\end{align}
%%%%%%%%%%%%%%%%%%%%%%%%%%%%%%%%%%%%%%%%%%%%%%%%%%%%%%%%%%%%%%%%%%%%%%%%%%%%%%%%%%%%%%%%%%%%
We now estimate each of the two terms on the right hand-side of the above equality. Combining arguments from Theorem \ref{thm4}, and adapting arguments from the proof of Proposition 3.8 (see estimates (64)--(68)) from \cite{anil7}, we have the following estimate, for the first term on the right hand-side of equation \eqref{prop7-eqn1}
\begin{align}\label{prop7-eqn2}
\sup_{z\in \gamma_{31}(\partial U_{k,a})}\big(1-|z|^{2}\big)^{k(n+1)+2}\Bigg|\sum_{\gamma\in \G}\frac{\ckn k(n+1)a_{ii}}{\big(Cz+D-\sum_{l=1}^{n}\big(\big(\sum_{m=1}^{n}a_{lm}z_{m}\big)+b_{l}\big)\overline{z_{l}}\big)^{k(n+1)}}\Bigg|\leq \notag\\[0.12cm]\sup_{z\in \gamma_{31}(\partial U_{k,a})}\sum_{\gamma\in \G}\frac{\ckn 8k(n+1)}{\cosh^{k(n+1)}\big(\dhyp(z,\gamma z)\slash 2 \big)}=O_{X_{\G}}\big(k^{n+3/2} \big).
\end{align} 
%%%%%%%%%%%%%%%%%%%%%%%%%%%%%%%%%%%%%%%%%%%%%%%%%%%%%%%%%%%%%%%%%%%%%%%%%%%%%%%%%%%%%%%%%%%%
Similarly, combining arguments from Theorem \ref{thm4}, and adapting arguments from the proof of Proposition 3.8 (see estimates (69) and (70)) from \cite{anil7}, we have the following estimate, for the second term on the right hand-side of equation \eqref{prop7-eqn1}
\begin{align}\label{prop7-eqn3}
&\sup_{z\in \gamma_{31}(\partial U_{k,a})}\big(1-|z|^{2}\big)^{k(n+1)+2}\times\notag\\&\Bigg|\sum_{\gamma\in \G}\frac{\ckn k(n+1)(k(n+1)+1)\big(\big(\sum_{l=1}^{n}a_{il}z_{l}\big) +b_{i}\big)\big( c_{i}-\big(\sum_{l=1}^{n}a_{li}\overline{z_{li}} \big)\big)}{\big(Cz+D-\sum_{l=1}^{n}\big(\big(\sum_{m=1}^{n}a_{lm}z_{m}\big)+b_{l}\big)\overline{z_{l}}\big)^{k(n+1)+2}}\Bigg|\leq \notag\\[0.12cm]&\sup_{z\in \gamma_{31}(\partial U_{k,a})}\sum_{\gamma\in \G}\frac{4\ckn k(n+1)(k(n+1)+1)}{\cosh^{k(n+1)}\big( \dhyp(z,\gamma z)\slash 2\big)}=O_{X_{\G}}\big(k^{n+5/2} \big).
\end{align}
%%%%%%%%%%%%%%%%%%%%%%%%%%%%%%%%%%%%%%%%%%%%%%%%%%%%%%%%%%%%%%%%%%%%%%%%%%%%%%%%%%%%%%%%%%%%
Finally, for any $\hat{z}:=(\hat{z}_1=\hat{x}_1+i\hat{y}_{1},\ldots,\hat{z}_n)^{t}\in \partial U_{k,a}$, using equation \eqref{gamma31}, we compute
\begin{align*}
1-\big|\gamma_{31}\hat{z}\big|^{2}=1-\bigg|\frac{\hat{z}_1+1}{\hat{z}_1-1}\bigg|-\sum_{j=2}^{n}\bigg|\frac{\hat{z}_{j}}{\hat{z}_1-1}\bigg|^{2}=\frac{-4\hat{x}_{1}-\sum_{j=2}^{n}|\hat{z}_{j}|^{2}}{\big((\hat{x}_1-1)^{2}+\hat{y}_{1}\big)^{2}},
\end{align*}
%%%%%%%%%%%%%%%%%%%%%%%%%%%%%%%%%%%%%%%%%%%%%%%%%%%%%%%%%%%%%%%%%%%%%%%%%%%%%%%%%%%%%%%%%%%%
which implies that
\begin{align}\label{prop7-eqn4}
\frac{1}{k}\ll_{X_{\G}}1-\big|\gamma_{31}\hat{z}\big|^{2}\ll_{X_{\G}}\frac{1}{k}.
\end{align}
%%%%%%%%%%%%%%%%%%%%%%%%%%%%%%%%%%%%%%%%%%%%%%%%%%%%%%%%%%%%%%%%%%%%%%%%%%%%%%%%%%%%%%%%%%%%
Combining estimates \eqref{prop7-eqn0}--\eqref{prop7-eqn4}, completes the proof of the proposition. 
\end{proof}
\end{prop}
%%%%%%%%%%%%%%%%%%%%%%%%%%%%%%%%%%%%%%%%%%%%%%%%%%%%%%%%%%%%%%%%%%%%%%%%%%%%%%%%%%%%%%%%%%%%

\vspace{0.1cm}
The following proposition, is an extension of Proposition 3.7 from \cite{anil7}.
%%%%%%%%%%%%%%%%%%%%%%%%%%%%%%%%%%%%%%%%%%%%%%%%%%%%%%%%%%%%%%%%%%%%%%%%%%%%%%%%%%%%%%%%%%%%

\vspace{0.1cm}
\begin{prop}\label{prop6}
With notation as above, for $k\gg1$, and $1\leq i \leq n$, and a fixed $a\in \mathbb{R}$ with $-2n\leq a\leq 2n$, , we have the following estimate
\begin{align}\label{prop6:eqn}
\sup_{z\in \overline{\calf}_{\G}}\big(1-|z|^{2}\big)^{k(n+1)-a}\bigg|\frac{\partial \bkn(z)}{\partial z_{i}} \bigg|=O_{X_{\G}}\big(k^{n+a+5/2} \big),
\end{align}
%%%%%%%%%%%%%%%%%%%%%%%%%%%%%%%%%%%%%%%%%%%%%%%%%%%%%%%%%%%%%%%%%%%%%%%%%%%%%%%%%%%%%%%%%%%%
where the implied constant depends only on $X_{\G}$.
%%%%%%%%%%%%%%%%%%%%%%%%%%%%%%%%%%%%%%%%%%%%%%%%%%%%%%%%%%%%%%%%%%%%%%%%%%%%%%%%%%%%%%%%%%%%
\begin{proof}
Let $\lbrace f_1,\ldots, f_{d_{k}}\rbrace$ denote an orthonormal basis of $\skn$. For any $z\in \calf_{\G}$, and $1\leq i,j \leq n$, and a fixed $a\in \mathbb{R}$ with $-2n\leq a\leq 2n$,   from the definition of the Bergman kernel, and Cauchy-Schwartz inequality, we have 
\begin{align}\label{prop8-eqn1}
\Bigg|\frac{\partial \bkn(z)}{\partial z_{i}}\Bigg| =\bigg|\sum_{l=1}^{d_{k}}\frac{\partial f_{l}(z)}{\partial z_{i}}\overline{f_{l}(z)}\bigg|\leq \sqrt{\sum_{l=1}^{d_{k}}\bigg|\frac{\partial f_{l}(z)}{\partial z_{i}}\bigg|^{2}}\cdot \sqrt{\sum_{l=1}^{d_{k}}\big|f_{l}(z)\big|^{2}}=\notag\\[0.12cm]
\sqrt{\Bigg|\frac{\partial^{2} \bkn(z)}{\partial z_{i}\partial \overline{z_{i}}}\Bigg|}\cdot\sqrt{\bkn(z)}.
\end{align}
%%%%%%%%%%%%%%%%%%%%%%%%%%%%%%%%%%%%%%%%%%%%%%%%%%%%%%%%%%%%%%%%%%%%%%%%%%%%%%%%%%%%%%%%%%%%
Combining Proposition \ref{prop8} with Theorem \ref{thm4}, we arrive at the following estimate
\begin{align}\label{prop8-eqn2}
&\sup_{z\in \overline{\calf}_{\G}}\big(1-|z|^{2}\big)^{k(n+1)-a}\Bigg|\frac{\partial \bkn(z)}{\partial z_{i}}\Bigg|\leq\notag\\[0.12cm]&\sup_{z\in \overline{\calf}_{\G}}\sqrt{\big(1-|z|^{2}\big)^{k(n+1)-2a}\Bigg|\frac{\partial^{2} \bkn(z)}{\partial z_{i}\partial \overline{z_{i}}}\Bigg|}\cdot \sup_{z\in \overline{\calf}_{\G}}\sqrt{\big|\bkn(z)\big|_{\mathrm{pet}}}=O_{X_{\G}}\big(k^{n+a+5/2} \big).
\end{align}
%%%%%%%%%%%%%%%%%%%%%%%%%%%%%%%%%%%%%%%%%%%%%%%%%%%%%%%%%%%%%%%%%%%%%%%%%%%%%%%%%%%%%%%%%%%%
Combining estimates \eqref{prop8-eqn1} and \eqref{prop8-eqn2}, completes the proof of the proposition. 
\end{proof}
\end{prop}
%%%%%%%%%%%%%%%%%%%%%%%%%%%%%%%%%%%%%%%%%%%%%%%%%%%%%%%%%%%%%%%%%%%%%%%%%%%%%%%%%%%%%%%%%%%%

\vspace{0.1cm}
The following proposition, is an extension of Proposition 3.7 from \cite{anil7}.
%%%%%%%%%%%%%%%%%%%%%%%%%%%%%%%%%%%%%%%%%%%%%%%%%%%%%%%%%%%%%%%%%%%%%%%%%%%%%%%%%%%%%%%%%%%%

\vspace{0.1cm}
\begin{prop}\label{prop7}
With notation as above, for $k\gg1$, and $1\leq i \leq n$, and a fixed $a\in \mathbb{R}$ with $-2n\leq a\leq 2n$, we have the following estimate
\begin{align}\label{prop7:eqn}
\sup_{z\in \overline{\calf}_{\G}}\big(1-|z|^{2}\big)^{k(n+1)-a}\bigg|\frac{\partial \bkn(z)}{\partial \overline{z}_{i}} \bigg|=O_{X_{\G}}\big(k^{n+a+5/2} \big),
\end{align}
%%%%%%%%%%%%%%%%%%%%%%%%%%%%%%%%%%%%%%%%%%%%%%%%%%%%%%%%%%%%%%%%%%%%%%%%%%%%%%%%%%%%%%%%%%%%
where the implied constant depends only on $X_{\G}$.
%%%%%%%%%%%%%%%%%%%%%%%%%%%%%%%%%%%%%%%%%%%%%%%%%%%%%%%%%%%%%%%%%%%%%%%%%%%%%%%%%%%%%%%%%%%%
\begin{proof}
The proof of the proposition follows directly from combining arguments from Proposition 3.7 from \cite{anil7} with Proposition \ref{prop6}.
\end{proof}
\end{prop}
%%%%%%%%%%%%%%%%%%%%%%%%%%%%%%%%%%%%%%%%%%%%%%%%%%%%%%%%%%%%%%%%%%%%%%%%%%%%%%%%%%%%%%%%%%%

\vspace{0.1cm}
Combining the above results, we now prove Main Theorem \ref{mainthm2}.
%%%%%%%%%%%%%%%%%%%%%%%%%%%%%%%%%%%%%%%%%%%%%%%%%%%%%%%%%%%%%%%%%%%%%%%%%%%%%%%%%%%%%%%%%%%

\vspace{0.1cm}
\begin{thm}\label{thm8}
With notation as above, for $k\gg1$, we have the following estimate
\begin{align}\label{thm8:eqn}
\sup_{z\in\overline{X}_{\G}}\bigg| \frac{\berkvol(z)}{\hypbnvol(z)}\bigg|=O_{X_{\G}}\big(k^{2(n-1)(n+2)+n+3}\big),
\end{align}
%%%%%%%%%%%%%%%%%%%%%%%%%%%%%%%%%%%%%%%%%%%%%%%%%%%%%%%%%%%%%%%%%%%%%%%%%%%%%%%%%%%%%%%%%%%
where the implied constant depends only on $X_{\G}$.
%%%%%%%%%%%%%%%%%%%%%%%%%%%%%%%%%%%%%%%%%%%%%%%%%%%%%%%%%%%%%%%%%%%%%%%%%%%%%%%%%%%%%%%%%%%
\begin{proof}
For any $z=(z_1,\ldots,z_n)^{t}\in\mathcal{F}_{\Gamma}$, from Lemma \ref{lem5}, and from the formula for the hyperbolic volume form, which as described in equation \eqref{defnhypmetric-vol}, we have
\begin{align}\label{thm8-eqn-1}
\bigg|\frac{\berkvol(z)}{\hypbnvol(z)}\bigg|=\big(1-|z|^{2} \big)^{n+1}\big|\mathrm{det}(M(z))\big|,
\end{align}
%%%%%%%%%%%%%%%%%%%%%%%%%%%%%%%%%%%%%%%%%%%%%%%%%%%%%%%%%%%%%%%%%%%%%%%%%%%%%%%%%%%%%%%%%%%
and the entires of the matrix 
\begin{align*}
M(z)=(m_{ij}(z))_{1\leq i, j \leq n},
\end{align*}
%%%%%%%%%%%%%%%%%%%%%%%%%%%%%%%%%%%%%%%%%%%%%%%%%%%%%%%%%%%%%%%%%%%%%%%%%%%%%%%%%%%%%%%%%%%
are as described in equations \eqref{lem5:eqn1} and \eqref{lem5:eqn2}, which we now estimate.
%%%%%%%%%%%%%%%%%%%%%%%%%%%%%%%%%%%%%%%%%%%%%%%%%%%%%%%%%%%%%%%%%%%%%%%%%%%%%%%%%%%%%%%%%%%

\vspace{0.1cm}
Using equations \eqref{lem5:eqn1} and \eqref{lem5:eqn2}, expanding along the first row of the matrix $M(z)$, we observe that
\begin{align}\label{thm8-eqn1}
&\big(1-|z|^{2} \big)^{n+1}\big|\mathrm{det}(M(z))\big|\ll_{X_{\G}} k^{n}+ \sum_{l=0}^{n-1}\frac{k^{l}g_{l}(z)}{\big( 1-|z|^{2}\big)^{2l-n-1}},
\end{align}
%%%%%%%%%%%%%%%%%%%%%%%%%%%%%%%%%%%%%%%%%%%%%%%%%%%%%%%%%%%%%%%%%%%%%%%%%%%%%%%%%%%%%%%%%%%
where
\begin{align*}
g_{l}(z):=\sum_{\mathcal{I}}\prod_{(\alpha,\beta)\in\mathcal{I}} \bigg(\frac{1}{\bkn(z)}\frac{\partial^{2}\bkn(z)}{\partial z_{\alpha}\partial\overline{z}_{\beta}}-
\frac{1}{\big(\bkn(z)\big)^{2}}\frac{\partial\bkn(z)}{\partial z_{\alpha}}\frac{\partial\bkn(z)}{\partial\overline{z}_{\beta}}\bigg),
\end{align*}
%%%%%%%%%%%%%%%%%%%%%%%%%%%%%%%%%%%%%%%%%%%%%%%%%%%%%%%%%%%%%%%%%%%%%%%%%%%%%%%%%%%%%%%%%%%
and $\mathcal{I}\subset  \big\lbrace (i,j)|\,1\leq i,j\leq n\rbrace$, and $|\mathcal{I}|$, the cardinality of the set $\mathcal{I}$ is equal to $n-l$.
%%%%%%%%%%%%%%%%%%%%%%%%%%%%%%%%%%%%%%%%%%%%%%%%%%%%%%%%%%%%%%%%%%%%%%%%%%%%%%%%%%%%%%%%%%%

\vspace{0.1cm}
For any $a\in\mathbb{R}$ with $-n\leq a\leq n$, from Propositions \ref{prop9} (we consider estimate \eqref{prop9:eqn2}, as the contribution of estimate \eqref{prop9:eqn2} dominates the contribution of estimate \eqref{prop9:eqn1}), and using the definition of map $\gamma_{13}$, for any $z\in\mathcal{F}_{\Gamma}$ we arrive at the following estimate
\begin{align*}
&\frac{1}{\big(1-|z|^{2}\big)^{a}}\bigg|\frac{1}{\bkn(z)}\frac{\partial^{2}\bkn(z)}{\partial z_{\alpha}\partial\overline{z}_{\beta}}\bigg|=\notag\\[0.12cm]&\frac{\big(1-|z|^{2}\big)^{k(n+1)-a}}{\big|\bkn(\gamma_{13}z)\big|_{\mathrm{pet}}}\cdot \bigg|\frac{\partial^{2}\bkn(z)}{\partial z_{\alpha}\partial\overline{z}_{\beta}}\bigg|\ll_{X_{\Gamma}} 
\frac{k^{-n/2}\big(1-|z|^{2}\big)^{k(n+1)-a}}{\big|\alpha(\gamma_{13}z)\big|^{(n-1)/2}}\cdot \bigg|\frac{\partial^{2}\bkn(z)}{\partial z_{\alpha}\partial\overline{z}_{\beta}}\bigg|,
\end{align*}
%%%%%%%%%%%%%%%%%%%%%%%%%%%%%%%%%%%%%%%%%%%%%%%%%%%%%%%%%%%%%%%%%%%%%%%%%%%%%%%%%%%%%%%%%%%
where 
\begin{align*}
\alpha(\gamma_{13}z)=\frac{2(1-|z|^{2})}{|z_{1}-1|^{2}}.
\end{align*}
%%%%%%%%%%%%%%%%%%%%%%%%%%%%%%%%%%%%%%%%%%%%%%%%%%%%%%%%%%%%%%%%%%%%%%%%%%%%%%%%%%%%%%%%%%%
Combining the above estimate with Proposition \ref{prop8}, we arrive at the following esimate
\begin{align}\label{thm8-eqn2}
&\frac{1}{\big(1-|z|^{2}\big)^{a}}\bigg|\frac{1}{\bkn(z)}\frac{\partial^{2}\bkn(z)}{\partial z_{\alpha}\partial\overline{z}_{\beta}}\bigg|\ll_{X_{\Gamma}}\notag\\[0.12cm]&k^{-n/2}\big(1-|z|^{2}\big)^{k(n+1)-(n-1)/2-a} \bigg|\frac{\partial^{2}\bkn(z)}{\partial z_{\alpha}\partial\overline{z}_{\beta}}\bigg|=O_{X_{\G}}\big( k^{n+a+4}\big).
\end{align}
%%%%%%%%%%%%%%%%%%%%%%%%%%%%%%%%%%%%%%%%%%%%%%%%%%%%%%%%%%%%%%%%%%%%%%%%%%%%%%%%%%%%%%%%%%%
Similarly, combining Proposition \ref{prop9} with Propositions  \ref{prop6} and \ref{prop7}, we have the following estimates
\begin{align*}
&\frac{1}{\big(1-|z|^{2}\big)^{a}}\bigg|\frac{1}{\bkn(z)}\bigg|\cdot\bigg|\frac{\partial\bkn(z)}{\partial z_{i_{\alpha}}}\bigg|=\notag\\[0.12cm]&\frac{\big(1-|z|^{2}\big)^{k(n+1)-a}}{\big|\bkn(z)\big|_{\mathrm{pet}}}\cdot \bigg|\frac{\partial\bkn(z)}{\partial z_{\alpha}}\bigg|=O_{X_{\G}}\big( k^{n+a+2}\big);
\end{align*}
%%%%%%%%%%%%%%%%%%%%%%%%%%%%%%%%%%%%%%%%%%%%%%%%%%%%%%%%%%%%%%%%%%%%%%%%%%%%%%%%%%%%%%%%%%%
and
\begin{align*}
&\frac{1}{\big(1-|z|^{2}\big)^{a}}\bigg|\frac{1}{\bkn(z)}\bigg|\cdot\bigg|\frac{\partial\bkn(z)}{\partial \overline{z}_{\beta}}\bigg|=\notag\\[0.12cm]&\frac{\big(1-|z|^{2}\big)^{k(n+1)-a}}{\big|\bkn(z)\big|_{\mathrm{pet}}}\cdot \bigg|\frac{\partial\bkn(z)}{\partial\overline{ z}_{\beta}}\bigg|=O_{X_{\G}}\big( k^{n+a+2}\big), 
\end{align*} 
%%%%%%%%%%%%%%%%%%%%%%%%%%%%%%%%%%%%%%%%%%%%%%%%%%%%%%%%%%%%%%%%%%%%%%%%%%%%%%%%%%%%%%%%%%%
using which we derive the following estimate
\begin{align}\label{thm8-eqn3}
\frac{1}{\big(1-|z|^{2}\big)^{a}}\bigg|\frac{1}{\big(\bkn(z)\big)^{2}}\frac{\partial\bkn(z)}{\partial z_{\alpha}}\frac{\partial\bkn(z)}{\partial\overline{z}_{\beta}}\bigg|=O_{X_{\G}}\big(k^{2n+a+ 4}\big).
\end{align}
%%%%%%%%%%%%%%%%%%%%%%%%%%%%%%%%%%%%%%%%%%%%%%%%%%%%%%%%%%%%%%%%%%%%%%%%%%%%%%%%%%%%%%%%%%%
Using estimates \eqref{thm8-eqn2} and \eqref{thm8-eqn3}, for any $0\leq l \leq n-1$, we derive
\begin{align}\label{thm8-eqn4}
&\frac{k^{l}g_{l}(z)}{\big( 1-|z|^{2}\big)^{2l-n-1}}\ll_{X_{\G}}\frac{k^{l}}{\big( 1-|z|^{2}\big)^{2l-n-1}}\times\notag\\[0.12cm]&\sum_{\mathcal{I}}\prod_{(\alpha,\beta)\in\mathcal{I}} \bigg(\frac{1}{\bkn(z)}\frac{\partial^{2}\bkn(z)}{\partial z_{\alpha}\partial\overline{z}_{\beta}}-\frac{1}{\big(\bkn(z)\big)^{2}}\frac{\partial\bkn(z)}{\partial z_{\alpha}}\frac{\partial\bkn(z)}{\partial\overline{z}_{\beta}}\bigg)\ll_{X_{\G}}\notag\\[0.12cm]&k^{l}\cdot k^{2n+4+2l-n-1}\cdot k^{(2n+4)(n-l-1)}=k^{2(n-1)(n+2)-l(2n+1)+n+3}.
 \end{align}
%%%%%%%%%%%%%%%%%%%%%%%%%%%%%%%%%%%%%%%%%%%%%%%%%%%%%%%%%%%%%%%%%%%%%%%%%%%%%%%%%%%%%%%%%%%
\vspace{0.1cm}
Combining estimates \eqref{thm8-eqn1} and \eqref{thm8-eqn4}, completes the proof of the theorem.
\end{proof}
\end{thm}
%%%%%%%%%%%%%%%%%%%%%%%%%%%%%%%%%%%%%%%%%%%%%%%%%%%%%%%%%%%%%%%%%%%%%%%%%%%%%%%%%%%%%%%%%%%
%%%%%%%%%%%%%%%%%%%%%%%%%%%%%%%%%%%%%%%%%%%%%%%%%%%%%%%%%%%%%%%%%%%%%%%%%%%%%%%%%%%%%%%%%%%

\vspace{0.2cm}
\subsection*{Acknowledgements} 
%%%%%%%%%%%%%%%%%%%%%%%%%%%%%%%%%%%%%%%%%%%%%%%%%%%%%%%%%%%%%%%%%%%%%%%%%%%%%%%%%%%%%%%%%%%
Both the authors thank A Raghuram for suggesting the problem. The first author is grateful to S. Das, J. Kramer, and A. Mandal for helpful discussions on estimates of automorphic forms. 
%%%%%%%%%%%%%%%%%%%%%%%%%%%%%%%%%%%%%%%%%%%%%%%%%%%%%%%%%%%%%%%%%%%%%%%%%%%%%%%%%%%%%%%%%%%
%%%%%%%%%%%%%%%%%%%%%%%%%%%%%%%%%%%%%%%%%%%%%%%%%%%%%%%%%%%%%%%%%%%%%%%%%%%%%%%%%%%%%%%%%%%
%%%%%%%%%%%%%%%%%%%%%%%%%%%%%%%%%%%%%%%%%%%%%%%%%%%%%%%%%%%%%%%%%%%%%%%%%%%%%%%%%%%%%%%%%%%
%%%%%%%%%%%%%%%%%%%%%%%%%%%%%%%%%%%%%%%%%%%%%%%%%%%%%%%%%%%%%%%%%%%%%%%%%%%%%%%%%%%%%%%%%%%

%%%%%%%%%%%%%%%%%%%%%%%%%%%%%%%%%%%%%%%%%%%%%%%%%%%%%%%%%%%%%%%%%%%%%%%%%%%%%%%%%%%%%%%%%%%
%%%%%%%%%%%%%%%%%%%%%%%%%%%%%%%%%%%%%%%%%%%%%%%%%%%%%%%%%%%%%%%%%%%%%%%%%%%%%%%%%%%%%%%%%%%
%%%%%%%%%%%%%%%%%%%%%%%%%%%%%%%%%%%%%%%%%%%%%%%%%%%%%%%%%%%%%%%%%%%%%%%%%%%%%%%%%%%%%%%%%%%
%%%%%%%%%%%%%%%%%%%%%%%%%%%%%%%%%%%%%%%%%%%%%%%%%%%%%%%%%%%%%%%%%%%%%%%%%%%%%%%%%%%%%%%%%%%
\end{document}